\documentclass{elsart}

\usepackage{amsmath,amssymb}
\usepackage{graphicx}
\usepackage{float}

\renewcommand{\P}{\mathbb{P}}
\newcommand{\E}{\mathbb{E}}
\newcommand{\R}{\mathbb{R}}
\newcommand{\Cov}{\mathrm{Cov}}
\newcommand{\ind}{{1\negthinspace \!\mathrm{I}}}
\newcommand{\eps}{\varepsilon}
\newcommand{\norm}[1]{\left\|#1\right\|}
\newtheorem{lemma}{Lemma}
\newtheorem{proposition}{Proposition}
\newtheorem{defi}{Definition}
\newenvironment{demo}{{ \em Proof~: }}{\hfill $\Box$ \smallskip}

\begin{document}

\begin{frontmatter}
% Title, authors and addresses

% use the thanksref command within \title, \author or \address for footnotes;
% use the corauthref command within \author for corresponding author footnotes;
% use the ead command for the email address,
% and the form \ead[url] for the home page:
% \title{Title\thanksref{label1}}
% \thanks[label1]{}
% \author{Name\corauthref{cor1}\thanksref{label2}}
% \ead{email address}
% \ead[url]{home page}
% \thanks[label2]{}
% \corauth[cor1]{}
% \address{Address\thanksref{label3}}
% \thanks[label3]{}

\title{Significant edges in the case of a non-stationary Gaussian noise}

\author[CEA]{I.~Abraham},
\ead{isabelle.abraham@cea.fr}
\author[Orleans]{R.~Abraham}, 
\ead{romain.abraham@univ-orleans.fr} 
\author[Paris5]{A.~Desolneux},
\ead{agnes.desolneux@math-info.univ-paris5.fr}
%\corauth[cor]{Corresponding author}
\author[Cachan]{S.~Li-Thiao-Te}
\ead{lithiaote@cmla.ens-cachan.fr}

\address[CEA]{CEA/DIF, 91680 Bruy\`eres le Chatel, France}
\address[Orleans]{Laboratoire MAPMO, F\'ed\'eration Denis Poisson, Universit\'e d'Orl\'eans, \\ B.P.~6759,
45067 Orl\'eans cedex 2, France} 
\address[Paris5]{Laboratoire MAP5, Universit\'e Ren\'e Descartes,
  45 rue des Saints-P\`eres, \\ 75270 Paris cedex 06, France}
\address[Cachan]{Laboratoire CMLA, ENS Cachan, 61 avenue du Pr\'esident Wilson, \\94235
  Cachan cedex, France}

% use optional labels to link authors explicitly to addresses:
% \author[label1,label2]{}
% \address[label1]{}
% \address[label2]{}

\author{}

\address{}

\begin{abstract}
% Text of abstract
In this paper, we propose an edge detection technique based on some local
smoothing of the image followed by a statistical hypothesis testing on the
gradient. 
An edge point being defined as a zero-crossing of the Laplacian, it is said to
be a significant edge point if the gradient at this point is larger than a
threshold $s(\eps)$ defined by: if the image $I$ is pure noise, then
$\P(\norm{\nabla I}\geq s(\eps) \bigm| \Delta I = 0) \leq\eps$. In other
words, a significant edge is an edge which has a very low probability to be
there because of noise. 
We will show that the threshold $s(\eps)$ can be explicitly computed in the case
of a stationary Gaussian noise. 
In images we are interested in, which are obtained by tomographic
reconstruction from a radiograph, this method fails since the Gaussian noise
is not stationary anymore. But in this case again, we will be able to give the
law of the gradient conditionally on the zero-crossing of the Laplacian, and
thus compute the threshold $s(\eps)$. We will end this paper with some
experiments and compare the results with the ones obtained with some other
methods of edge detection.    
\end{abstract}

\end{frontmatter}

\section{Introduction}

This work is part of some specific physical experiments which consist
in studying radially symmetric objects \cite{Cho}. These objects are composed of
several materials and one of the interesting features is the location
of the frontier between the different materials.

\begin{figure}[H]
\begin{center}
\includegraphics[width=5cm]{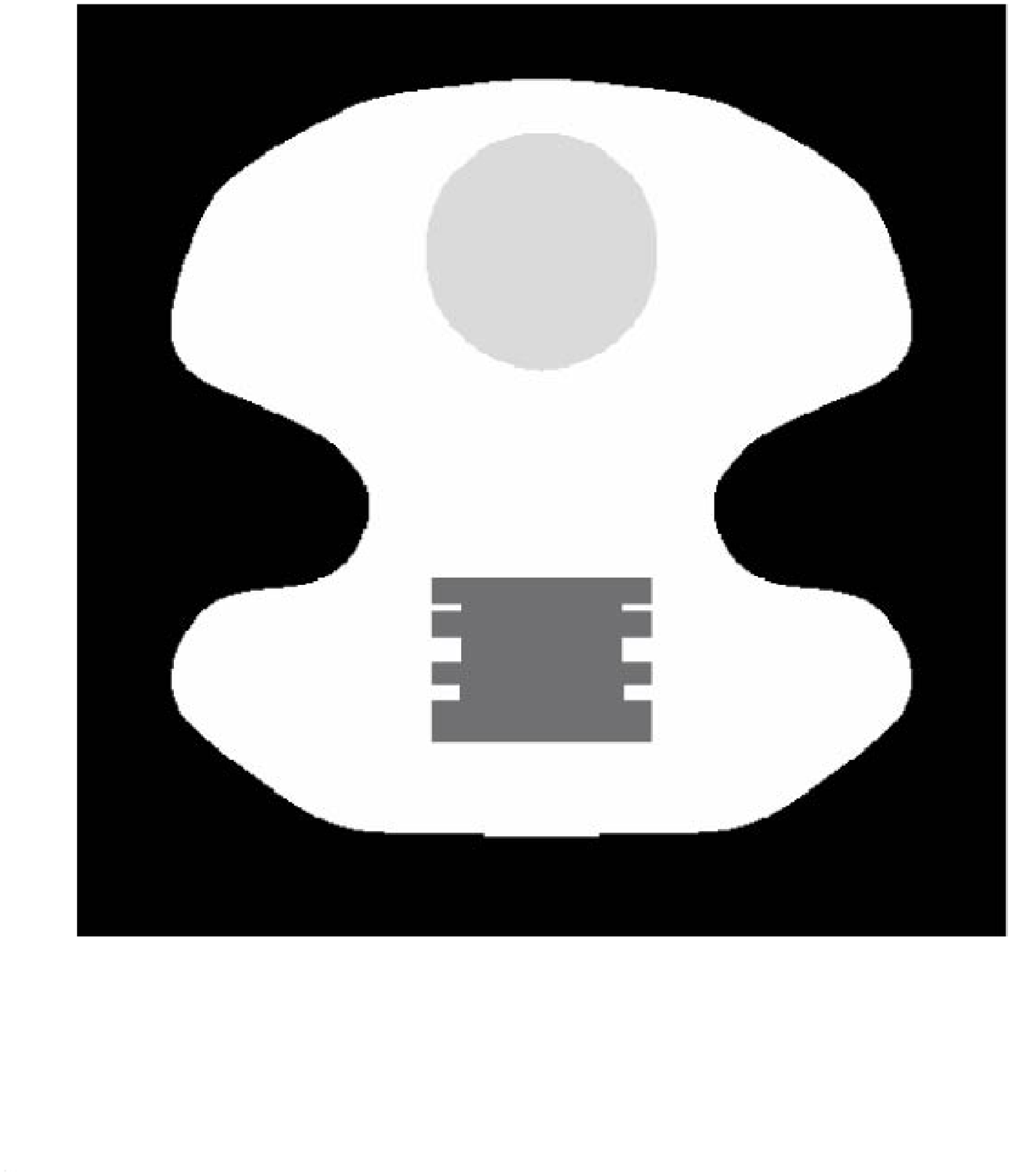}
\end{center}
\caption{\label{fig:object}Slice of a studied object.}
\end{figure}

To
describe such an object, it is enough to give the densities of the
materials on a slice of the object that contains the
symmetry axis. An example of studied object is given on Figure
\ref{fig:object}.

\begin{figure}[H]
\begin{center}
\begin{tabular}{ccc}
\includegraphics[width=4.2cm]{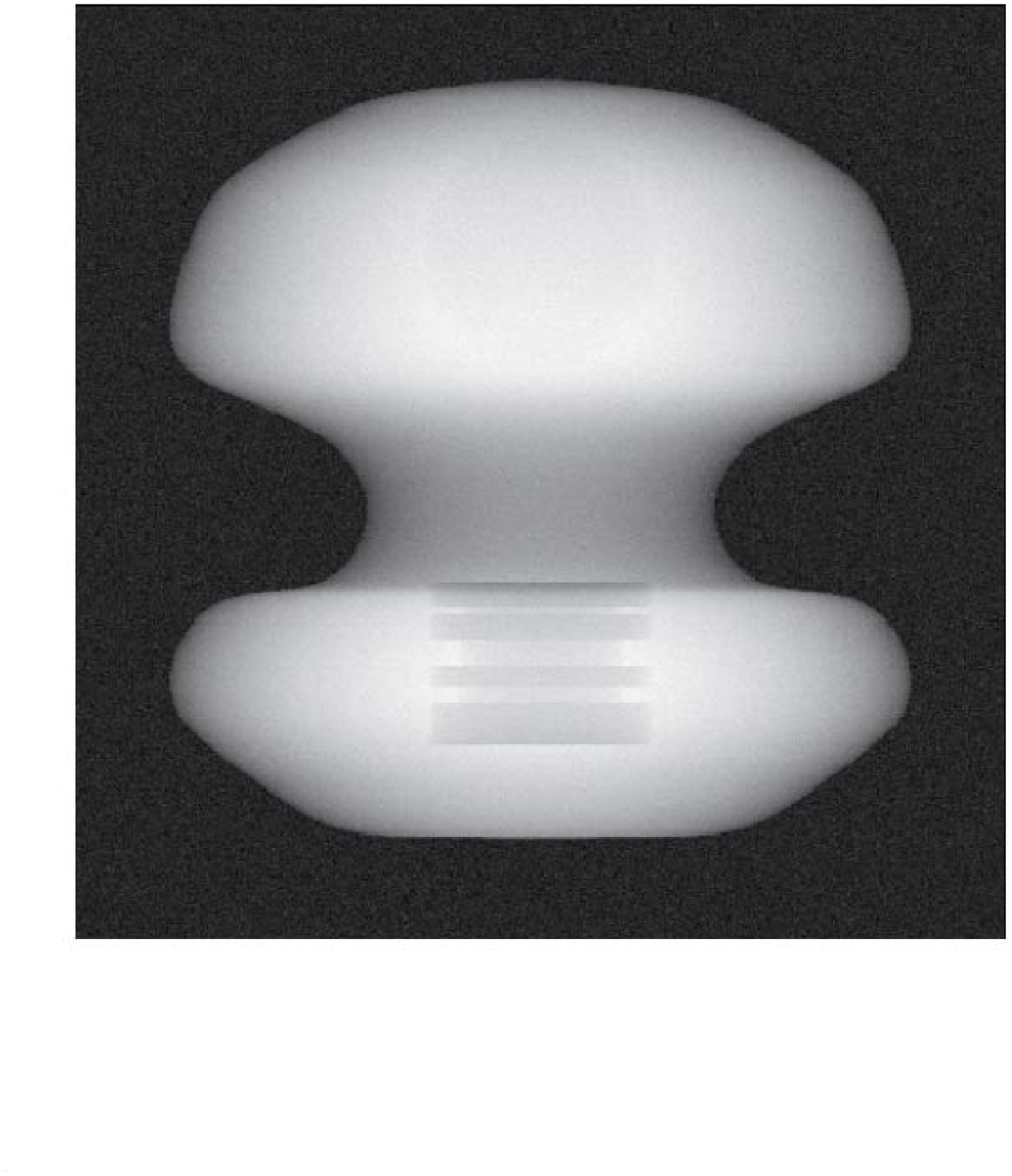} &
\includegraphics[width=4.2cm]{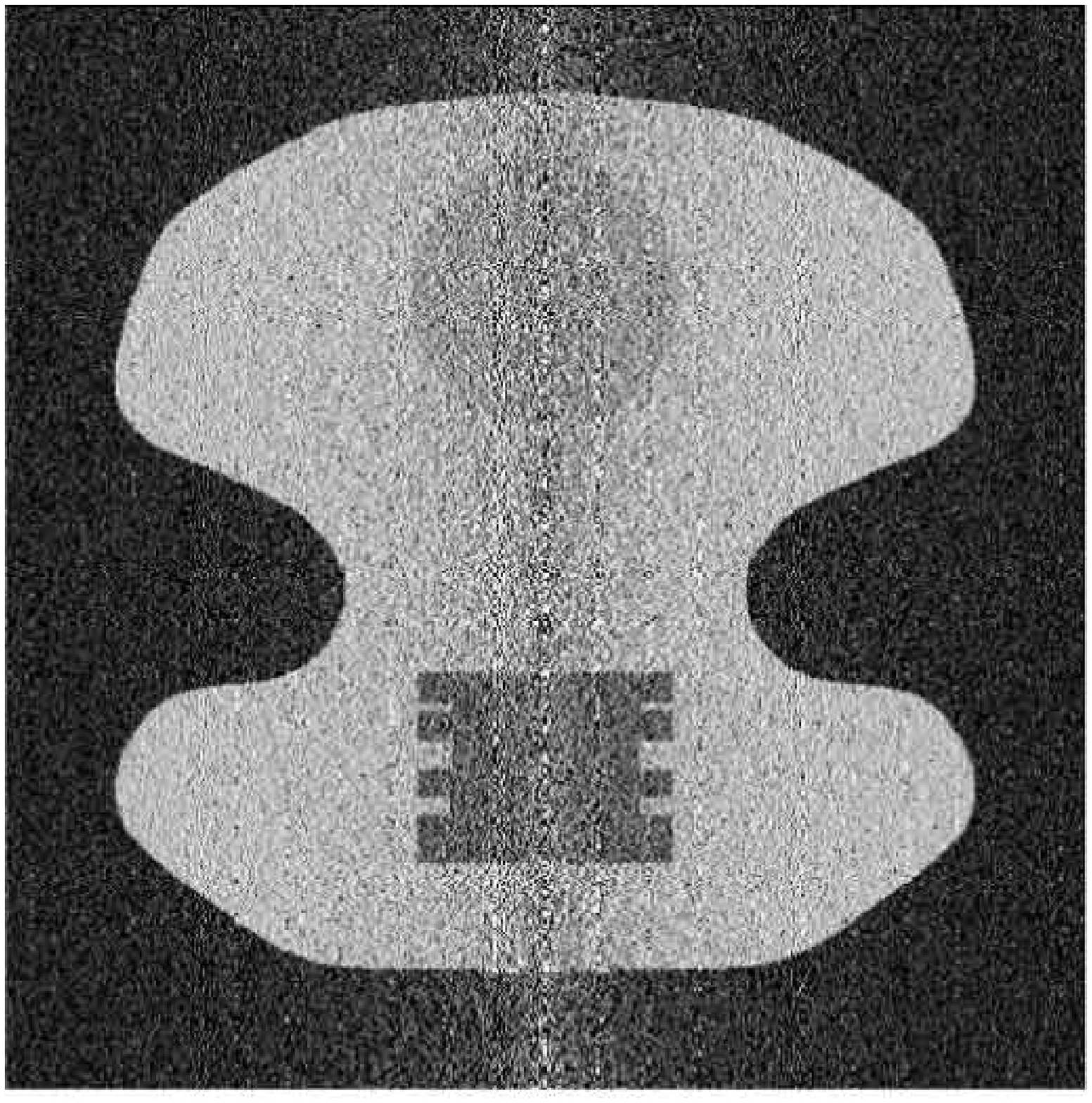} &
\includegraphics[width=4.2cm]{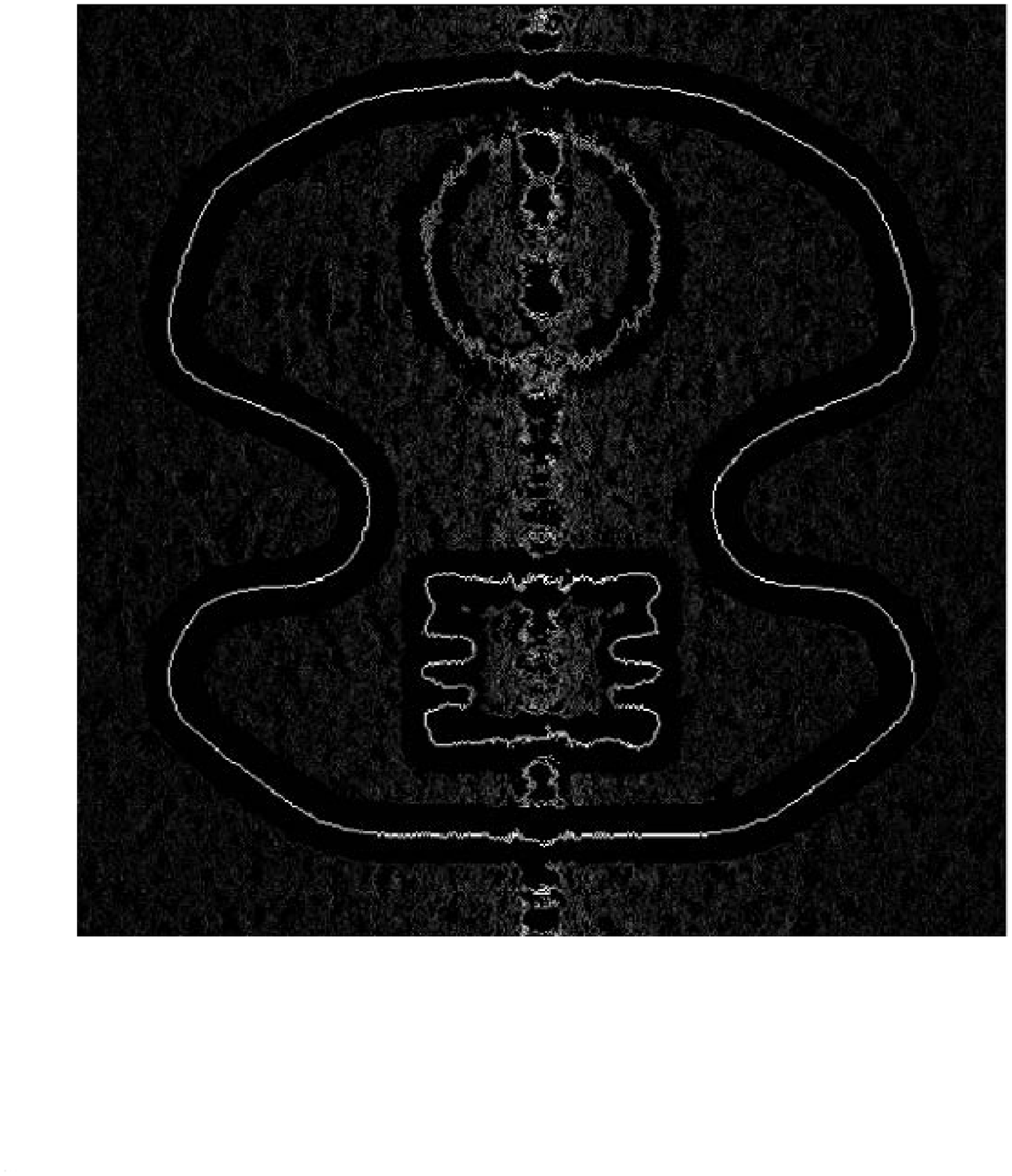}\\
(a) & (b) & (c)
\end{tabular}
\end{center}
\caption{\label{fig:treatment}
(a)  Radiograph of the object of Figure \ref{fig:object}.
(b)  Tomographic reconstruction.
(c)  Edge detection on the tomography.
}
\end{figure}

To look at the interior of this object, a radiography is performed
(see Figure \ref{fig:treatment}(a)), then a tomographic reconstruction is
computed (Figure \ref{fig:treatment}(b)) and finally an edge detection is
made (Figure \ref{fig:treatment}(c)).

At this point, let us mention that the tomographic reconstruction is
not an usual one. Indeed, the usual inverse Radon transform (and the
usual reconstruction algorithms such as filtered back-projection)
operates on a slice of the object that is orthogonal to the rotation
axis. Here, because of  the radial symmetry assumption, the
reconstruction can be global \cite{Din}. This reconstruction will be detailed
in Section \ref{tomography.subsec}.

As we can see on Figure \ref{fig:treatment}(c), many detected
edges do not correspond to real features. This is due to the high
level of noise.
For the time being, the selection of the edges is manually executed. The
goal of this work is to perform this selection automatically. For that
purpose, the edge detector will not be changed but we will compute
also other significant features that will allow us to select the
``true'' edges. 

The ideas used here mainly come from previous work of Desolneux,
Moisan and Morel \cite{DMM}. Informally speaking, they define the
notion of significant edges by computing the probability that some
edge-related events appear in an image of pure noise. When this
probability is small enough, the edge is probably a feature of the
image and not due to the noise. Unfortunately, their method assumes that the 
noise is stationary which, as easily seen on Image 
\ref{fig:treatment}(b), is not the case in our study because of the tomographic
inversion (see Section \ref{expconstant.subsec} for some examples of results
obtained with their method). Moreover, their study is quite general and, apart
from the stationarity, no assumption is made on the noise. 

In our case, as we deal with specific images, the noise is well-known
and some statistical models can be used. Indeed, we may suppose that
the noise on the radiograph (\ref{fig:treatment}(a)) is a Gaussian
white noise with mean zero and with a variance that can easily be
estimated. Then, a tomographic inversion is performed. As this operation
is linear, we still obtain a Gaussian noise but it is now correlated
and non-stationary.

The edge detector will not be modified here. It consists in estimating
the Laplacian at each point, and edge points are then defined as the
zero-crossings of the Laplacian. As we already said, we only add some
features that will discriminate the significant edges. The easiest
feature to compute is a kind of contrast measurement $C$ based on a
gradient estimate.
To be more precise, we consider an image $I$ of pure noise (that is a
realization of our model of noise after tomographic reconstruction), estimate
the gradient and the Laplacian of $I$ at a 
point $(u,v)$ (with an abuse of notation, we will denote by
$\nabla I(u,v)$ and $\Delta I(u,v)$ these estimates and by $C(u,v)$
the contrast value) and we compute,
for a fixed $\varepsilon>0$, the smallest value $s(\eps)$ for which
\begin{equation}\label{eq:proba}
\P(C(u,v)\ge s(\eps)\bigm| \Delta I(u,v)=0)\le
\varepsilon.
\end{equation}
Then, we perform an edge detection on the studied image $f$ (where we
also estimate $\nabla f$ and $\Delta f$ by the same method) and we keep the
points $(u,v)$ of the studied image $f$ that 
satisfy
\begin{itemize}
\item $\Delta f(u,v)=0$ (an edge is present at point $(u,v)$).
\item $C(u,v)\ge s(\eps)$ (this edge is significant).
\end{itemize}
>From a statistical point of view, this consists in performing an
hypothesis test. We consider a point $(u,v)$ where an edge takes place
($\Delta (u,v)=0$) and we test the null hypothesis ``the zero-crossing
of the Laplacian is due to the noise''.
The level of the test $\varepsilon$ is arbitrarily chosen and
related to the number of false detections allowed. It will be set to 
$\eps=10^{-5}$ hereafter. Let us mention that the threshold value $s(\eps)$
varies slowly with respect to $\varepsilon$. For instance, in the case
of a white noise (see Section 3), the threshold value can be computed
explicitly and is proportional to $\sqrt{-\ln\varepsilon}$. When the null
hypothesis is rejected, the edge is retained as it comes from a
``true'' feature of the image, whereas when the null hypothesis is
accepted, the zero-crossing of the Laplacian may come from the noise
and the edge is not meaningful.

Let us mention at this point that such statistical approaches
have already been used for edge detection in \cite{TLB}, \cite{QB} or
\cite{MR}. They usually use estimates of the gradient
based on finite differences which fail in our case. Moreover, the
noise is in most cases stationary. Let us also cite \cite{CMS} where the
authors have 
modified the method of \cite{DMM} to take into account the
non-stationarity of some images, by a local noise estimate. Their
work is still general and does not make any assumption on the noise
structure. As we deal with specific experiments, the noise is always
the same and well-known and we can take proper advantage of this knowledge.

The paper is organized as follows: in Section 2, we present the edge
detector based on the estimate of the gradient and the
Laplacian. Then, in Section 3, our method is presented in the case of a
Gaussian white noise. Of course, this does not correspond to our case
but the computations are easier and show the performance of this
method. In Section 4, we will first describe the tomographic inversion
and the operators involved, and then describe the noise model we have
to deal with. We will then apply the significant edges
detection method in the framework of this non-stationary noise. We will end
the section with some experiments and comparisons with other methods.

\section{Estimating the Gradient and  the Laplacian}
\label{sec:estimations}

In this section, we introduce a method for edge detection. We consider
that the image is a real-valued function $(u,v)\mapsto f(u,v)$ of two
continuous real  
parameters $u$ and $v$. Then, we say that there exists an edge at point
 $(u,v)$ if the  
Laplacian of $f$ is zero at this point. Moreover, the computation of
the contrast function $C$ will be based on the gradient of $f$ (see
the end of this section for the choice of this function). As the images are
very noisy, these derivatives cannot be 
estimated by usual finite differences. The method used here, sometimes
known as Savitsky-Golay smoothing, consists in locally approximating
the image by a polynomial. The derivatives of the polynomial are then
identified with those of the image.

\subsection{An optimization problem}

Let $(u,v)$ denote the point where we want to compute the first and
second order derivatives of the image $f$. We choose 2 parameters : $d$
which is the maximum degree of the approximating polynomial and $r$
which is the radius of the ball on which we perform the
approximation. We denote by $B_r(u,v)$ the ball of radius $r$ centered
at point $(u,v)$. We will simply write $B_r$ when the center of the ball is the
origin $(0,0)$ of $\mathbb{R}^2$. We are then looking for a polynomial $P$ of
degree less that $d$ such that 
\begin{equation}
E(P)=\int_{B_r}\bigl( f(u+x,v+y)-P(x,y)\bigr)^2\,dx\,dy
\label{energy.eq}
\end{equation}
is minimal among all polynomials of degree less than $d$.
In other words, we are looking for the best approximation of $f$ by a
polynomial  of degree less than $d$ on the ball
$B_r(u,v)$ in the sense of the $L^2$-norm. \\
This is an optimization problem where the unknowns are the coefficients of
the polynomial. As the problem is convex, there is a unique solution
(given by the orthogonal projection of $f$ on the space of polynomials
of degree less than $d$) which is easily computed by solving the
equations
$$\frac{\partial E}{\partial a_i}=0$$
where the $a_i$'s denote the coefficients of the polynomial.

{\it Role of the ball radius.}
Two parameters are arbitrary chosen in this method. The first one is
the ball radius $r$. The larger $r$ is, the more effective the
smoothing is.
 The influence of the noise is therefore attenuated with a
large $r$ but the location of the edge is then less precise. We
must consequently make a balance between noise smoothing and edge
detection accuracy. For instance, if we have a small level of noise or
if the edges are very complicated (with high curvature), we must
choose a small value for $r$.

{\it Role of the polynomial degree.}
The second parameter is the polynomial degree. Here again a large
value of $d$ gives a better approximation but does not smooth the noise
enough. In fact, as we are, in a first step, interested in the points where
the Laplacian is zero, it appears that a second-order polynomial is
enough. Of course, the estimate of the first order derivatives with a
polynomial of degree $2$ is not very good and highly depends on the
size of the window $B_r$. But we will see that this drawback can be useful for the choice of a contrast function.

In what follows, the approximation is made with a polynomial of degree $d=2$,
and the first and second order derivatives of the image are
identified with those of the approximating polynomial.

\subsection{Computations with a second order polynomial}

Let us first introduce some notations. In the following, we will set
$$\forall i,j \in \mathbb{N} , \hspace{0.3cm} b_{ij}(r)= \int_{B_r}x^iy^jdx\, dy.$$
As the ball $B_r$ is symmetric, we have that $b_{ij}(r)=0$ as soon as
$i$ or $j$ is odd and that $b_{ij}(r)=b_{ji}(r)$ for all $i$,$j$.
In order to have simple expressions, we will also set:
$$b(r) = b_{20}(r) , \hspace{0.2cm} \alpha (r) =-\frac{2b_{20}(r)}{b_{00}(r)} \hspace{0.2cm} \text{ and } \hspace{0.2cm}
\beta (r) =\frac{1}{2}\left(b_{40}(r)+b_{22}(r)-\frac{2b_{20}^2(r)}{b_{00}(r)}\right) $$

\begin{lemma}
The gradient and the Laplacian of the polynomial of degree $2$ which is the
best approximation of $f$ on the ball $B_r(u,v)$ for the $L^2$-norm, being
respectively denoted by  
$\nabla_r f(u,v) = ( \frac{\partial_r f}{\partial x}(u,v) , \frac{\partial_r
  f}{\partial y}(u,v) )$  and $\Delta_r f(u,v)$, are given by: 
\begin{align*}
\frac{\partial_r f}{\partial x}(u,v) & =\frac{1}{b(r)}\int_{B_r}x\,
f(u+x,v+y)\, dx\, dy\\ 
\frac{\partial_r f}{\partial y}(u,v) & =\frac{1}{b(r)}\int_{B_r}y\,
f(u+x,v+y)\, dx\, dy\\ 
\Delta_r f(u,v) & = \frac{1}{\beta(r)}
\int_{B_r}f(u-x,v-y)\left(\alpha(r)+x^2+y^2\right)\, dx\, dy.
\end{align*}
\end{lemma}

\begin{demo}

We consider a polynomial of degree $2$ which we write
$$P(x,y)=a_{xx}x^2+a_{yy}y^2+a_{xy}xy+a_xx+a_yy+a_0.$$
The equations obtained by writing $\nabla E(P)=0$, where $E(P)$ is given by Equation (\ref{energy.eq}), are: 
$$\left\{\begin{array}{lll}
b_{40}(r)a_{xx}+b_{22}(r)a_{yy}+b_{20}(r)a_0 & = &
\displaystyle\int_{B_r}x^2\,f(u+x,v+y)\,dx\,dy\\
b_{22}(r)a_{xx}+b_{40}(r)a_{yy}+b_{20}(r)a_0 & = & \displaystyle\int_{B_r}y^2\,
f(u+x,v+y)\,dx \,dy\\
b_{22}(r)a_{xy} & = & \displaystyle\int_{B_r} x\, y\, f(u+x,y+v)\, dx\, dy\\
b_{20}(r)a_x & = & \displaystyle\int_{B_r}x\, f(u+x,y+v)\, dx\, dy\\
b_{20}(r)a_y & = & \displaystyle\int_{B_r}y\, f(u+x,v+y)\, dx\, dy\\
b_{20}(r)a_{xx}+b_{20}(r)a_{yy}+b_{00}(r)a_0 & = & \displaystyle\int_{B_r}f(u+x,v+y)\, dx\, dy
\end{array}\right.$$

We then obtain the following estimates for the derivatives:
\begin{align*}
\frac{\partial P}{\partial x}(0,0) & = a_x =\frac{1}{b_{20}(r)}\int_{B_r}x\, f(u+x,v+y)\, dx\, dy\\
\frac{\partial P}{\partial y}(0,0) & = a_y =\frac{1}{b_{20}(r)}\int_{B_r}y\, f(u+x,v+y)\, dx\, dy\\
\Delta P(0,0) & = 2(a_{xx}+a_{yy})  \\
& =\frac{2}{b_{40}+b_{22}-\frac{2b_{20}^2}{b_{00}}}
\int_{B_r}f(u+x,v+y)\left(-\frac{2b_{20}}{b_{00}}+x^2+y^2\right)\, 
dx\, dy.
\end{align*}
\end{demo}

\subsection{Choice of the contrast function}

We would like to use a contrast function based on the estimates of the first
and second derivatives of the image $f$ obtained in the previous section.

The simplest contrast function we can choose is the norm of the gradient:
$$C_1(u,v)=\norm{\nabla_r f(u,v)}.$$
Indeed, the value of this norm tells how sharp the edge is. This
contrast function is efficient and will be used when the images we
deal with are piecewise constant.

\begin{figure}[H]
\begin{center}
\includegraphics[width=5cm]{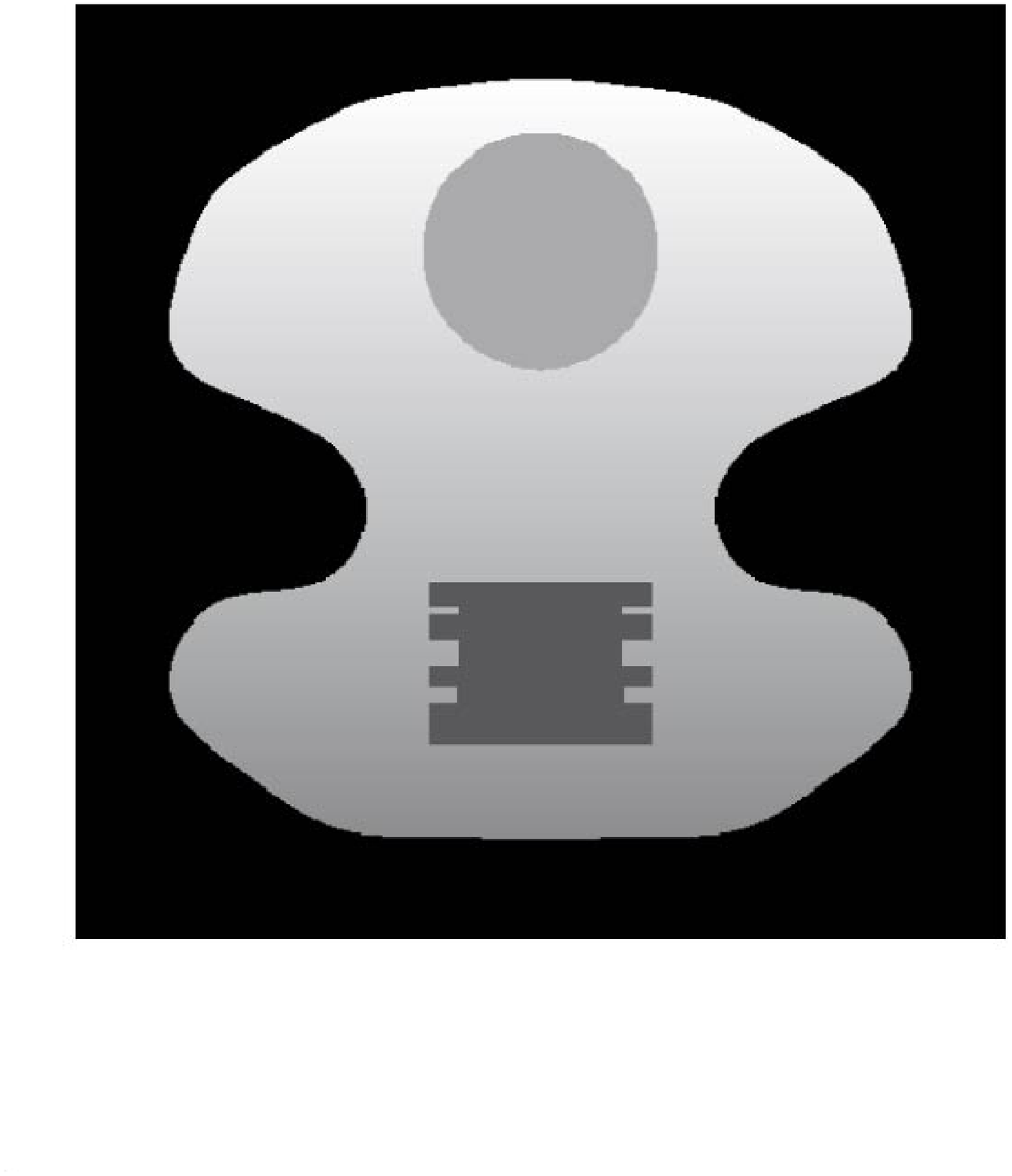}
\end{center} 
\caption{\label{fig:slope}Object with an inhomogeneous
  material}
\end{figure}

\begin{figure}[H]
\begin{center}
\includegraphics[width=5cm]{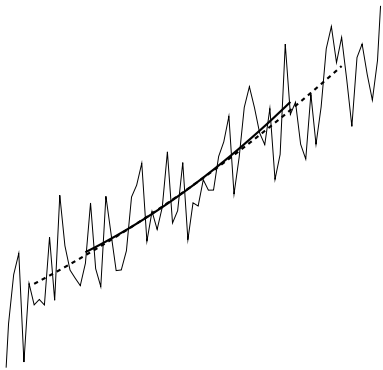}
\end{center}
\caption{\label{fig:approx_slope}A noisy constant slope: the gradient of the approximating
  polynomial does not depend on the radius $r$.}
\end{figure}

However, in many cases, the objects we handle are not homogeneous and
their images contain some slopes (see Figure~\ref{fig:slope}). In this
case, the gradient norm is not a good contrast function. Indeed, let
us consider an image with a constant slope with some noise (see Figure
\ref{fig:approx_slope}). We would 
like to say that no edge is significant in that case. However, the
value of the gradient norm (which will be close to the value of the
slope) will always be greater that the threshold value $s$ when the
noise level is small.

\begin{figure}[H]
\begin{center}
\includegraphics[width=5cm]{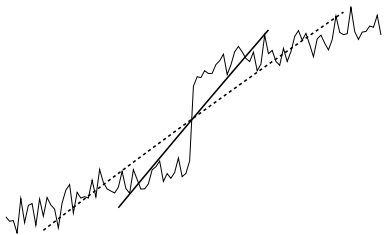}
\end{center}
\caption{\label{fig:approx_slope_edge}An edge on a noisy slope: the approximating
  polynomials with two different values of the radius $r$.}
\end{figure}

In the latter case, we take advantage of the dependence of the first
order derivatives estimates with respect to the ball radius. Indeed,
the estimates of the gradient in the case of the constant slope in
Figure \ref{fig:approx_slope} will not depend on the size of the window (see
Figure \ref{fig:approx_slope}) whereas, when an edge (a discontinuity)
occurs, the estimates do depend on that radius (see Figure
\ref{fig:approx_slope_edge}). So, we can use as a contrast function
the function
$$C_2(u,v)=\norm{\nabla_{r_1}f(u,v)-\nabla_{r_2}f(u,v)}$$
where $r_1<r_2$ and $\nabla_r f$ denotes the value of the gradient
estimate with a ball of radius $r$.

\section{Significant edges in the case of a Gaussian white noise}
\subsection{White noise and Wiener integral}
We recall here the definition and the main properties of a white noise
in a continuous setting and of the Wiener integral. We refer to
\cite{Wal}, \cite{Hid} or \cite{HKPS} for more on white noise and the
Wiener integral.

\begin{defi}
A Gaussian white noise on $\R^2$ of variance $\sigma^2$ is a random
function $W$ defined on the Borel sets $A$ of $\R^2$ of finite Lebesgue
measure (denoted by $|A|$) such that
\begin{itemize}
\item $W(A)$ is a Gaussian random variable (r.v.) with mean $0$ and variance $\sigma^2|A|$ ,
\item If $A_1\cap A_2=\emptyset$, the r.v. $W(A_1)$ and $W(A_2)$ are
  independent and
$$W(A_1\cup A_2)=W(A_1)+W(A_2).$$
\end{itemize}
\end{defi}

Such a function $W$ exists but is not a true measure since the
two-parameters process
$$B(s,t):=W\bigl((0,s]\times(0,t]\bigr)$$
(usually called the Brownian sheet) is of unbounded total variation.

Nevertheless we can define the so-called Wiener integral $\int fdW$
for every function $f$ in $L^2(\R_+^2)$. We can also define the
derivatives of the Brownian sheet in the sense of Schwartz
distributions (although the Brownian sheet is nowhere
differentiable). Thus, we define
$$\dot B(s,t)=\frac{\partial ^2B(s,t)}{\partial s\partial t}$$
and we have
$$\int fdW=\int _{\R_+^2}f(u,v)\dot B(u,v)du\, dv\qquad a.s.$$
for every function $f$ in the Schwartz space.

With a slight abuse of notations, we call $\dot B$ a Gaussian white
noise and we always denote by $\int _{\R_+^2}f(u,v)\dot B(u,v)du\, dv$
the Wiener integral with respect to this white noise, for every
function $f\in L^2$. The main properties of this integral are
\begin{itemize}
\item For every $f$, the r.v. $\displaystyle \int _{\R_+^2}f(u,v)\dot
  B(u,v)du\, dv$ is a Gaussian r.v. with mean 0 and variance
$\sigma^2\displaystyle \int _{\R_+^2}f(u,v)^2du\,dv.$
\item For every $f,g$, the random vector
$$\left(\int_{\R^2}f(u,v)\dot B(u,v)du\,dv,\int_{\R^2}g(u,v)\dot
  B(u,v)du\,dv\right)$$
is Gaussian with covariance
$$\sigma^2\int _{\R^2}f(u,v)g(u,v)du\,dv.$$
\end{itemize}

We will use these properties to compute the laws of $\nabla I$ and
$\Delta I$.

\subsection{Laws of the gradient and of the Laplacian}
We suppose here that our noise is a Gaussian white noise, of variance $\sigma^2$. As we have already said, this case is not the one
we are interested in and our method is probably over-performed by other
standard methods in that case. The goal of this section is to present
our method in a simple case where the computations are easy to do and can be
carried out in a continuous setting. We will only focus here on the
case of piecewise constant objects and therefore we will use the contrast
function $C_1$.

\begin{lemma}
If the image $I$ is a Gaussian white noise of variance $\sigma^2$, then 
$$\left(\frac{\partial_r I}{\partial x},\frac{\partial_r I}{\partial
y},\Delta_r I\right)$$
is a Gaussian vector with mean zero and covariance matrix
$$\left(\begin{array}{ccc}
\frac{\sigma^2}{b(r)} & 0 & 0\\
0 & \frac{\sigma^2}{b(r)} & 0\\
0 & 0 & V(r,\sigma)
\end{array}\right)  , \hspace{0.2cm} \text{ where } \hspace{0.2cm}  V(r,\sigma)=\frac{\sigma^2}{\beta^2(r)}\int_{B_r}\left(\alpha(r)+x^2+y^2\right)^2dx\, dy .$$
\label{bruitblanc.lem}
\end{lemma}

\begin{demo}
We compute the laws of the approximate derivatives of $I$ when $I=\dot B$. We
recall that these derivatives are given by 
\begin{align*}
\frac{\partial_r I}{\partial x}(u,v) & =\frac{1}{b(r)}\int_{B_r}x\,
\dot B(u+x,v+y)\, dx\, dy\\
\frac{\partial_r I}{\partial y}(u,v) & =\frac{1}{b(r)}\int_{B_r}y\,
\dot B(u+x,v+y)\, dx\, dy\\
\Delta_r I(u,v) &
=\frac{1}{\beta(r)}\int_{B_r}\dot B(u+x,v+y)\left(\alpha(r)+x^2+y^2\right)\,
dx\, dy.
\end{align*}
Because of the stationarity of $\dot B$, they have the same law as
\begin{align*}
\frac{\partial_r I}{\partial x}(0,0) & =\frac{1}{b(r)}\int_{B_r}x\,
\dot B(x,y)\, dx\, dy\\
\frac{\partial_r I}{\partial y}(0,0) & =\frac{1}{b(r)}\int_{B_r}y\,
\dot B(x,y)\, dx\, dy\\
\Delta_r I(0,0) &
=\frac{1}{\beta(r)}\int_{B_r}\dot B(x,y)\left(\alpha(r)+x^2+y^2\right)\,
dx\, dy.
\end{align*}

As we deal with Wiener integrals, we deduce that the vector
$$\left(\frac{\partial_r I}{\partial x},\frac{\partial_r I}{\partial
y},\Delta_r I\right)$$
is a Gaussian vector with mean zero.

To compute its covariance matrix,
let us recall that, if $X$ and $Y$ are random variables defined by
\begin{align*}
X & =\int_{B_r}h_1(x,y)\dot B(x,y)\, dx\, dy\\
Y & =\int_{B_r}h_2(x,y)\dot B(x,y)\, dx\, dy
\end{align*}
then we have
$$\Cov (X,Y)=\sigma^2\int_{B_r}h_1(x,y)h_2(x,y)\, dx\, dy.$$

Consequently, we have for instance:
$$\Cov\left(\frac{\partial_r I}{\partial x},\frac{\partial_r I}{\partial
y}\right)  = \frac{\sigma^2}{b^2(r)}\int _{B_r}x\, y\, dx\, dy 
= 0. $$
By some analogous calculations, we finally get the following covariance matrix
for our Gaussian vector: 
$$\left(\begin{array}{ccc}
\frac{\sigma^2}{b(r)} & 0 & 0\\
0 & \frac{\sigma^2}{b(r)} & 0\\
0 & 0 & V(r,\sigma)
\end{array}\right)$$
where
$$V(r,\sigma)=\frac{\sigma^2}{\beta^2(r)}\int_{B_r}\left(\alpha(r)+x^2+y^2\right)^2dx\,
dy.$$
\end{demo}

Thanks to this lemma, we immediately have the following properties:

\begin{itemize}
\item The random variable $\left\|\nabla_r I\right\|^2$ is the sum of
two squared independent Gaussian random variables which have the same variance. It is
therefore distributed as a $\chi^2$-law. More precisely, its law is
$$\frac{\sigma^2}{b(r)}\chi^2(2)$$
where $\chi^2(2)$ denotes a $\chi^2$-law with two degrees of freedom.
\item The random variable $\Delta_r I$ is a Gaussian random variable with mean
zero and variance $V(r,\sigma)$.
\item The random variables $\left\|\nabla_r I\right\|^2$ and $\Delta_r I$
are independent.
\end{itemize}

\subsection{Computation of the threshold}

\begin{proposition}
Let $I$ be a Gaussian white noise and
let $s(\eps)$ be the threshold value such that 
$$\P\bigl(\left\|\nabla_r I\right\| \ge s(\eps) \bigm| \Delta_r I=0\bigr)\le
\varepsilon.$$ 
Then
$$s(\eps)=\sqrt{-\frac{2\sigma^2}{b(r)}\ln\varepsilon}.$$
\end{proposition}

\begin{demo}
To begin with, as the random variables $\left\|\nabla_r I\right\|^2$ and $\Delta_r I$
are independent, we can forget the conditioning and only compute
$$\P\bigl(\left\|\nabla_r I\right\| \ge s(\eps)\bigr) =
\P\bigl(\left\|\nabla_r I\right\|^2\ge s(\eps)^2\bigr).$$ 
As a consequence of Lemma \ref{bruitblanc.lem}, we have that the law of
$\left\|\nabla_r I\right\|^2$ is $\frac{\sigma^2}{b(r)}\chi^2(2)$. 
Now, since the density of a $\chi^2(2)$ law is the one of a
$\Gamma\left(\frac{1}{2},1\right)$ law, we have that the law of
$\left\|\nabla_r I\right\|^2$ is given by 
$$\P\bigl(\left\|\nabla_r I\right\|^2 \ge s^2\bigr) =
\int_{\frac{b(r)}{\sigma^2}s^2}^{+\infty} \frac{1}{2}e^{-\frac{t}{2}}dt = \exp
\left( - \frac{b(r) s^2}{2\sigma^2} \right). 
 $$
This finally leads to the announced threshold value $s(\eps)$.
\end{demo}

\subsection{Experiments}\label{sec:experiments_white}

We consider the piecewise constant object of Figure \ref{fig:object}
with some additive Gaussian white noise. The densities of the different
materials of this object are: 
\begin{itemize}
\item $1$ for the outer material,
\item $0.8$ for the material inside the circle,
\item $0.3$ for the other inner material.
\end{itemize}

\begin{figure}[H]
\begin{center}
\begin{tabular}{ccc}
\includegraphics[width=4.2cm]{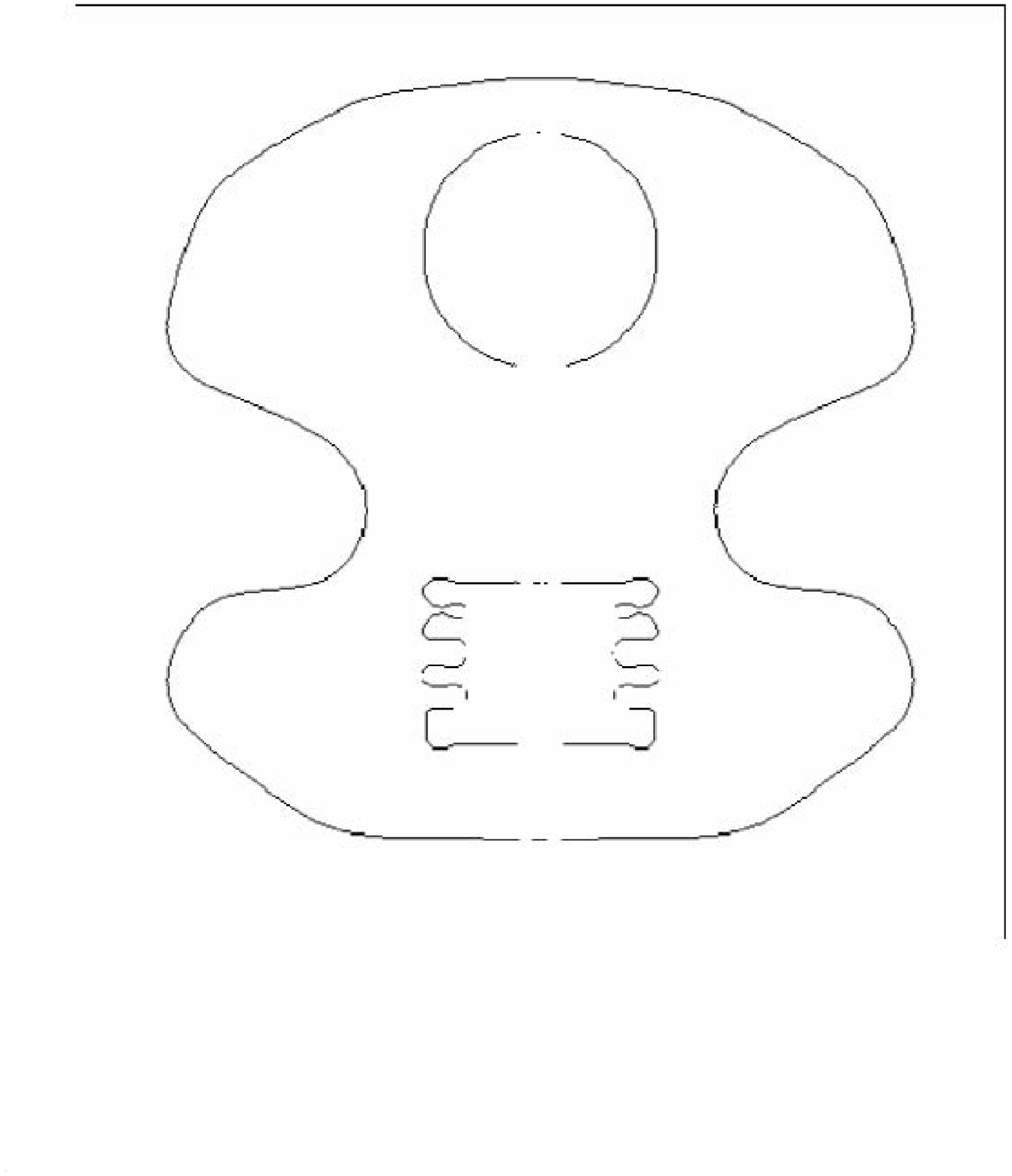} &
\includegraphics[width=4.2cm]{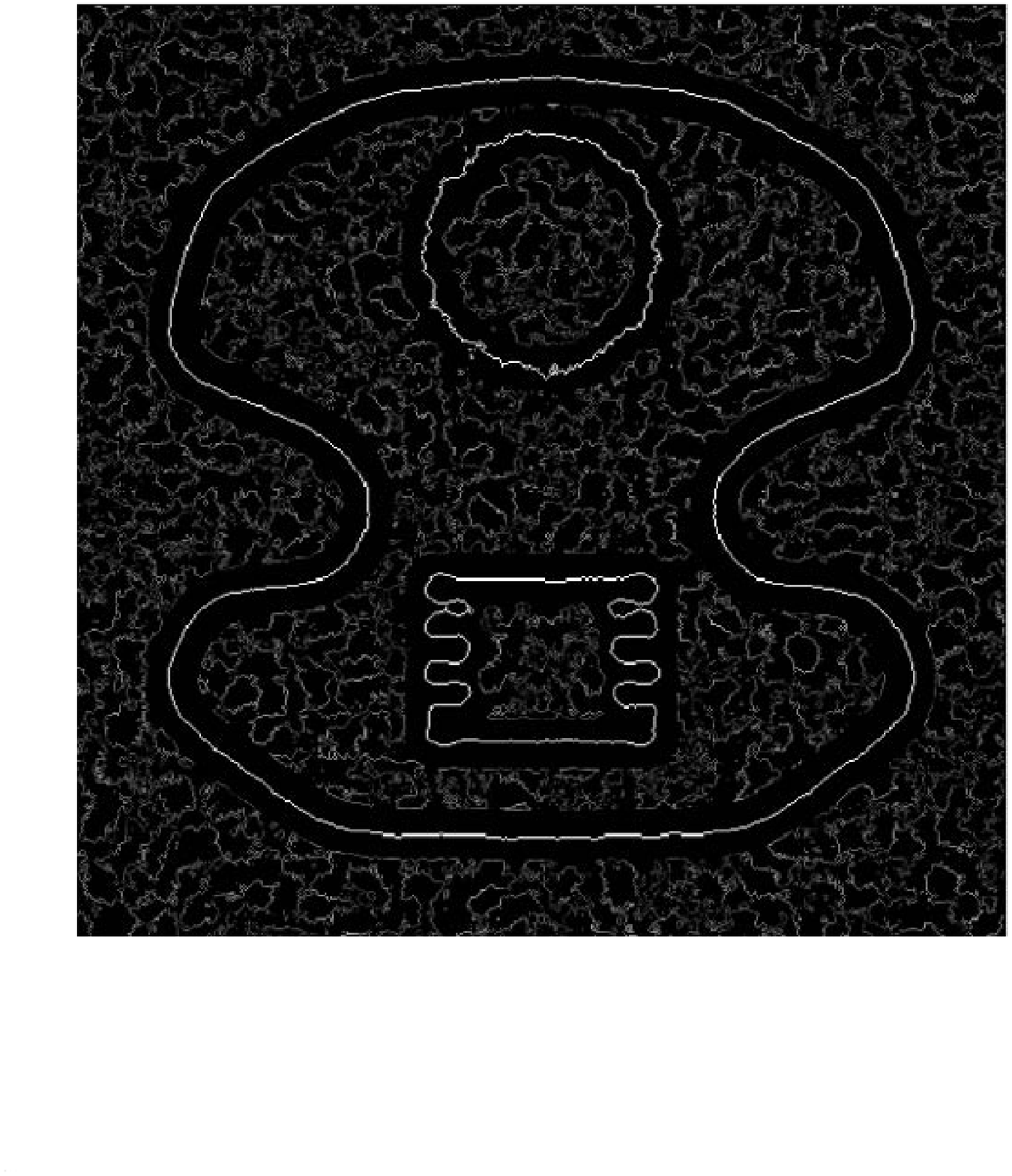} &
\includegraphics[width=4.2cm]{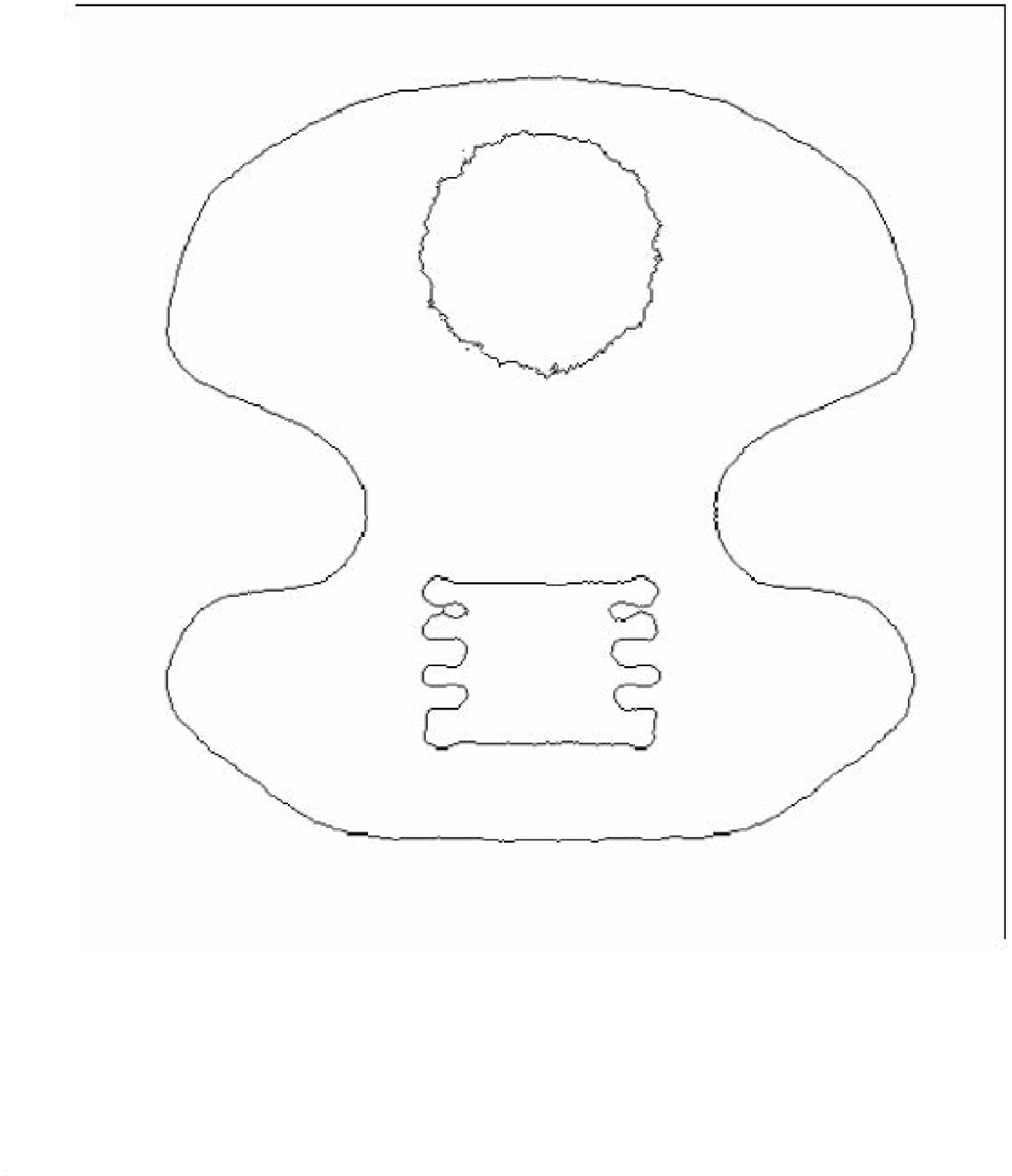}\\
(a) & (b) & (c)
\end{tabular}
\end{center}
\caption{\label{fig:experiments_white_1} (a) The noisy image ($\sigma=0.2$).
(b) The zero-crossings of the Laplacian with the contrast function
  $C_1$ visualized in grey-level.
(c) The extracted significant edges ($\eps=10^{-5}$).
}
\end{figure}

\begin{figure}[H]
\begin{center}
\begin{tabular}{ccc}
\includegraphics[width=4.2cm]{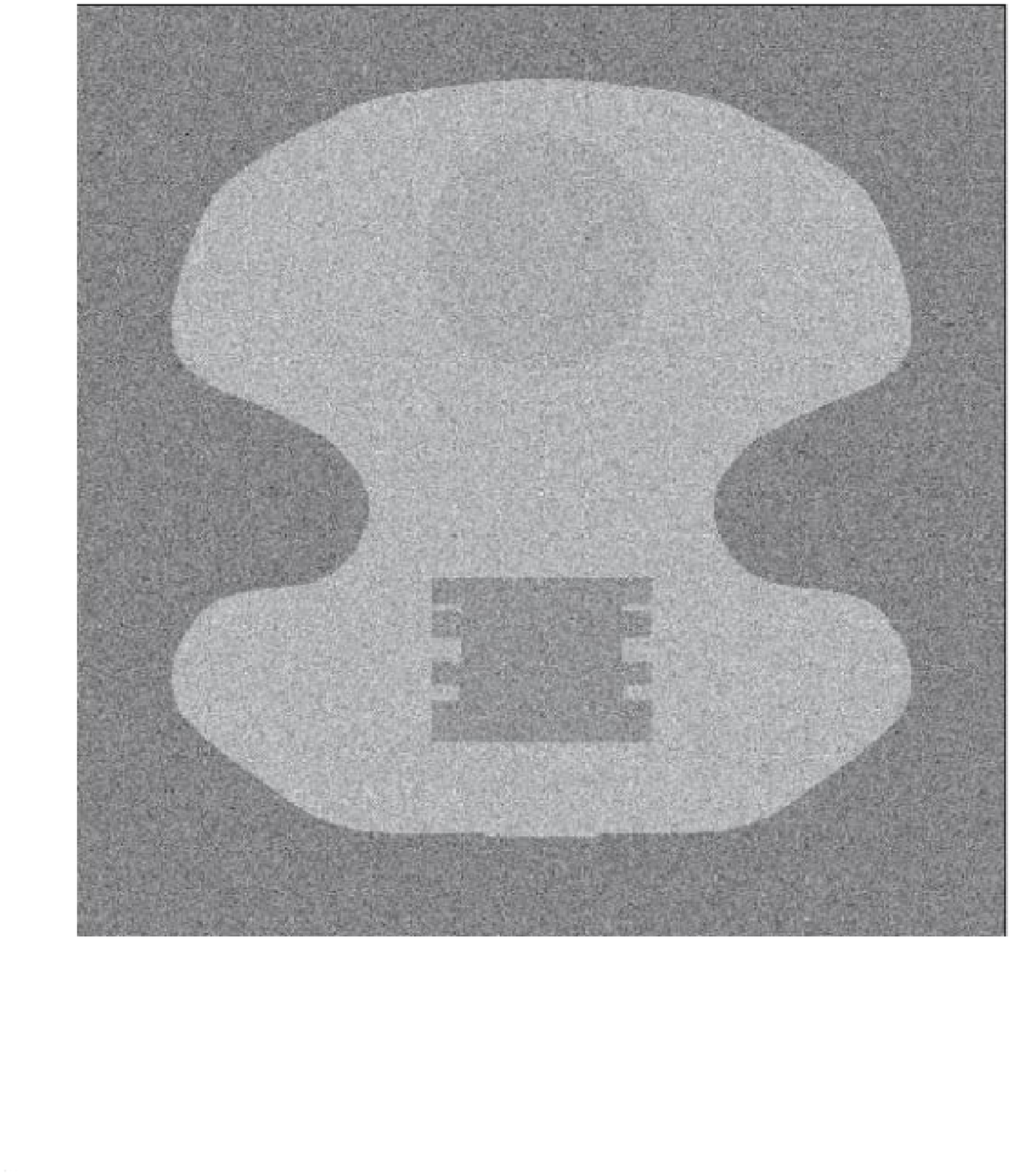} &
\includegraphics[width=4.2cm]{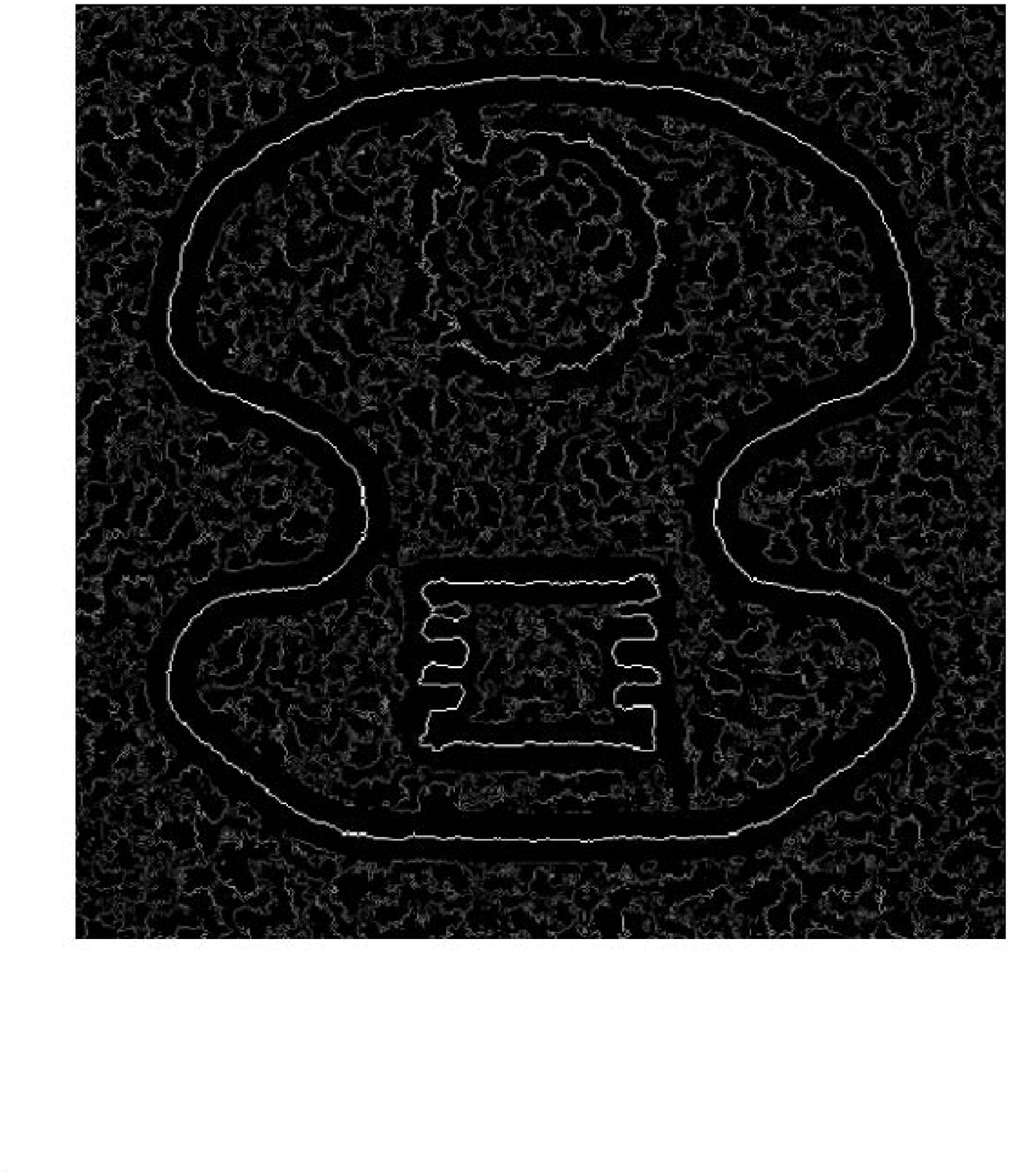} &
\includegraphics[width=4.2cm]{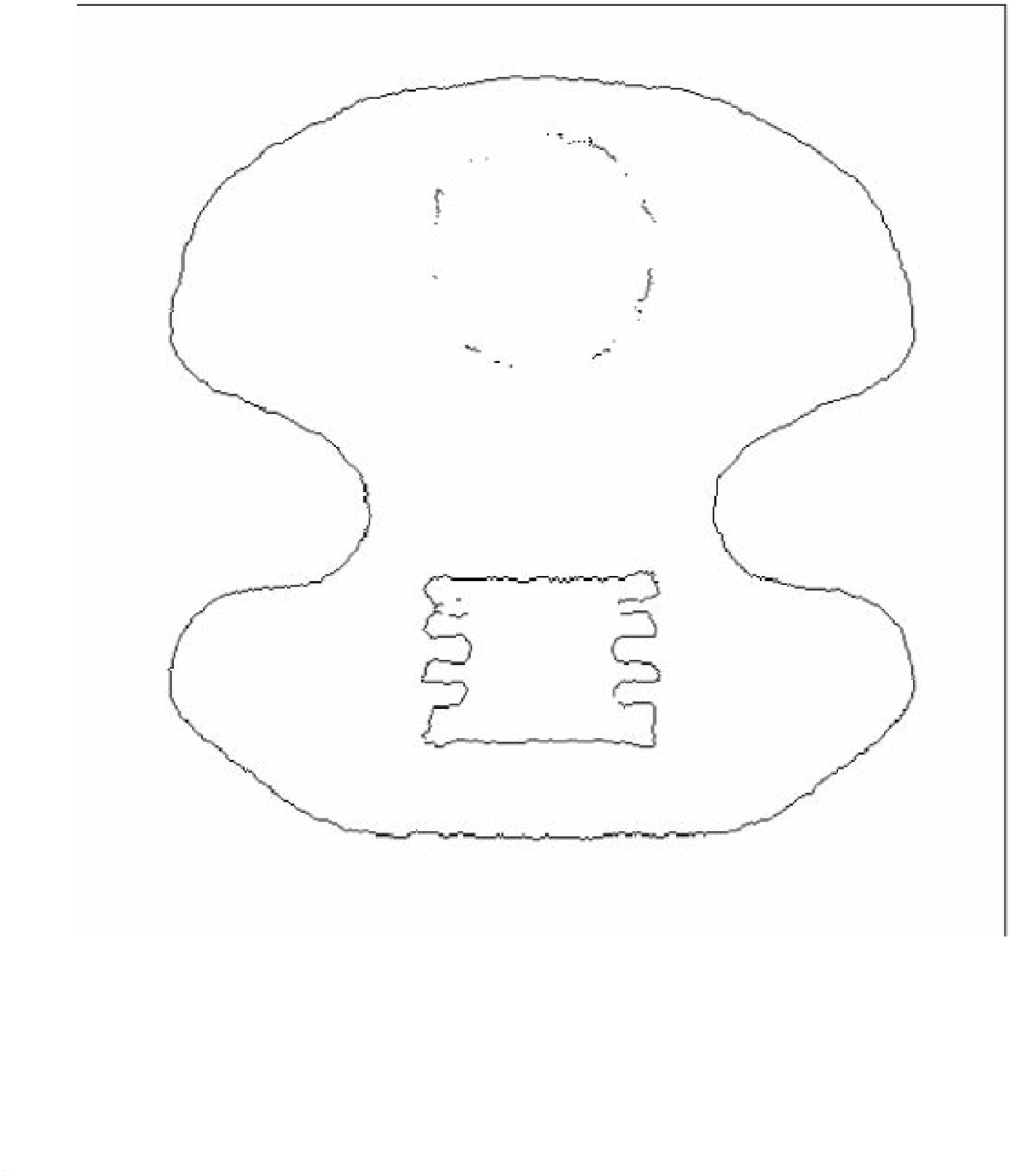}\\
(a) & (b) & (c)
\end{tabular}
\end{center}
\caption{\label{fig:experiments_white_2} (a) The noisy image ($\sigma=0.4$).
(b) The zero-crossings of the Laplacian with the contrast function
  $C_1$ visualized in grey-level. (c) The extracted significant edges ($\eps=10^{-5}$).
}
\end{figure}

The standard deviation of the Gaussian noise is $\sigma=0.2$ in the experiments
of Figure \ref{fig:experiments_white_1} and is $\sigma=0.4$ in the
experiments of Figure \ref{fig:experiments_white_2}. Both images have the same
size $512\times 512$ pixels. 
The experiments have been carried out with a ball of radius $r=12$ pixels.\\
The different images of Figures \ref{fig:experiments_white_1} and \ref{fig:experiments_white_2} are respectively:
\begin{itemize}
\item (a) The noisy image.
\item (b) The zero-crossings of the Laplacian with the contrast function
  $C_1$ visualized in grey-level (the white color corresponds to 
  high value for the contrast function $C_1$).
\item (c) The extracted significant edges ($\eps=10^{-5}$).
\end{itemize}

In the case of a signal-to-noise ratio large enough (Figure
\ref{fig:experiments_white_1}), all the edges are well detected and the
``false'' edges are removed. Let us nevertheless mention that, with our method,
the edges which have a high curvature are smoothed. This drawback is even more
important 
when the ball radius $r$ is large (the influence of the value of this
radius will be studied in the experiments of the next section).

When the noise level is rather large (Figure
\ref{fig:experiments_white_2}), some edges of the image cannot be extracted from
the noise (it happens when the contrast associated to this edge is close to the
noise level).

\section{Significant edges in the case of a Gaussian white noise on the
  radiograph}

\subsection{Tomography}
\label{tomography.subsec}

Let us turn now to the more realistic case we are interested in. As we
mentioned it in the introduction, we first make a radiography of an
object. Tomography is the inverse problem associated with
reconstructing the initial object from its radiograph. This is now a
well-known problem as it is the key tool in medical scanner imagery (or other
medical imaging systems).

To begin with, let us describe what a radiography is from a mathematical
point of view. The studied object is exposed to $X$-rays that go through
it. Some of the $X$-photons are absorbed. As an output, we observe the quantity
of $X$-photons that have not been absorbed by the material, and we thus measure
in some sense the ``mass'' of material the ray went through. More
precisely, if the object is described by its density $\mu$ (which is
a function of the space coordinates), what can be measured at some
point of the receptor is
$$\int_{ray}\mu\, d\ell$$
where ``ray'' means the straight line that goes from the source to the
studied point of the receptor (we suppose that the X-rays source is just a
point, which implies that the previous line is unique).

We also assume that the X-rays source is far away from the object so
that the rays are assumed to be parallel. Then, to reconstruct
any object from its radiographs, we must turn around the object and make a
radiography for every angle $\theta\in[0,\pi)$. This leads to the
so-called Radon transform of the object, which is known
to be invertible. This is the principle of the medical scanner.

\begin{figure}[H]
\begin{center}
\includegraphics[width=8cm]{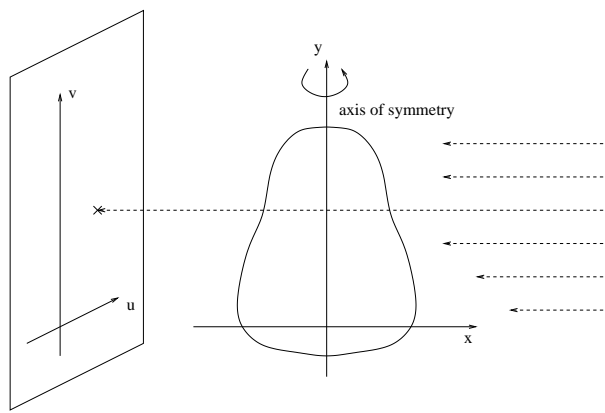}
\end{center}
\caption{\label{tomo.fig} Radiography of a radially symmetric object.}
\end{figure}

In our case, as the object is radially symmetric, if we turn around
the object with for rotation axis the symmetry axis of the object, all the
radiographs are exactly the same. Consequently, a single radiograph of such an
object 
is enough to perform the tomographic reconstruction. Indeed, if
$f(x,y)$ denotes the density along a slice that contains the symmetry
axis (see Figures \ref{fig:object} and \ref{tomo.fig}), then a radiograph of
this object is given by
$$g(u,v)=2\int_{|u|}^{+\infty}f(x,v)\frac{x}{\sqrt{x^2-u^2}}dx.$$
This is a linear transform and we will denote it hereafter by 
$$g=Hf.$$
As we already said, this linear operator $H$ is invertible and we in
fact know explicitly its inverse on the space of continuously
differentiable functions $g$:
$$f(x,y)=(H^{-1}g)(x,y)=-\frac{1}{\pi}\int_x^{+\infty}
\frac{1}{\sqrt{u^2-x^2}}\frac{\partial g}{\partial u}(u,y)du.$$

Our assumption on the noise is that it is an additive Gaussian white
noise on the radiograph (i.e. on $g$). But what we want is to study the object
given by $f$. So we must transform the white noise by the operator
$H^{-1}$. Unfortunately, because of the singularity of the integral at
$x=0$, we cannot apply the operator $H^{-1}$ to a white noise $\dot
B$, even in a $L^2$-sense. Therefore, we will work in a discrete
framework: the images $f$ and $g$ are naturally discretized (as they
are numerical images). This leads to a discretization of the operator
$H$, which we will still denote by $H$ and which now may be viewed as
a matrix. The discretization is made in such a way that the symmetry axis
($x=0$) is settled between two pixels so that the previous singularity
does not appear. 
This matrix is then invertible and we denote by $H^{-1}$ its
inverse which we can make now operate on a discrete Gaussian white noise.

\subsection{Law of the noise on the tomographic reconstruction}

Let us consequently consider a field $\eta=(\eta_{i,j})_{{1\le i\le
    p},{1\le j\le n}}$ of
i.i.d. random Gaussian variables with mean $0$ and variance $\sigma^2$. Let
us define $I=(I_{i,j})$ the random field obtained after tomographic
reconstruction i.e. after making $H^{-1}$ operate on
$\eta=(\eta_{i,j})$. In fact, as the X-rays are supposed to be
parallel, the reconstruction can be made line by line independently
and therefore, if we consider the row vectors
$$\vec\eta_i=(\eta_{i,1},\ldots ,\eta_{i,n}) \hspace{0.2cm} \text{ and }
\hspace{0.2cm} \vec I_i=(I_{i,1},\ldots,I_{i,n})$$
then, there exists an invertible matrix $M$ (independent of $i$, and of size $n\times n$)   such that
$$\vec I_i=\vec \eta_i M.$$

Consequently, the law of $I$ is characterized by the following
properties:
\begin{itemize}
\item $I=(I_{i,j})$ is a Gaussian random field.
\item For $i\ne k$, $\vec I_i$ and $\vec I_k$ are independent.
\item For each $i$, the vector $\vec I_i$ is a Gaussian vector of mean $0$ and
  covariance matrix
$$\Gamma=\sigma^2 M^tM ,$$
where $M^t$ denotes the transpose of $M$.
\end{itemize}

\subsection{Laws of the gradient and of the Laplacian}

The expressions obtained in Section \ref{sec:estimations} for the
gradient and for the Laplacian of an image in a continuous setting are
easily translated in the discrete framework we now deal with. Indeed,
we have
\begin{align*}
\frac{\partial_r I}{\partial x}(u,v) & =\frac{1}{b(r)}\sum_{(i,j)\in
  B_r}jI_{u+i,v+j}\\
\frac{\partial_r I}{\partial y}(u,v) & =\frac{1}{b(r)}\sum_{(i,j)\in
  B_r}iI_{u+i,v+j}\\
\Delta_r I(u,v) &
 = \frac{1}{\beta(r)}\sum_{(i,j)\in B_r} \left(\alpha(r)+i^2+j^2\right)I_{u+i,v+j}
\end{align*}
where $B_r$ now denotes the discrete ball of radius $r$ i.e.
$$B_r=\{(i,j),\ i^2+j^2\le r^2 \}$$
and where the constants $\alpha(r)$, $\beta(r)$, $b(r)$, $\ldots$ are the
discrete analogous of the constants of Section \ref{sec:estimations}.

With these estimates, the contrast functions $C_1$ and $C_2$ are
easily comptuted. They are both of the form
$$C(u,v)=\sqrt{C_x^2(u,v)+C_y^2(u,v)}$$
with
$$C_x(u,v)  =\sum_{i,j} j {c_{ij}} I_{u+i,v+j} \hspace{0.2cm} \text{ and } \hspace{0.2cm}
C_y(u,v) = \sum_{i,j} i {c_{ij}} I_{u+i,v+j} , $$
where the coefficients $c_{ij}$ are given by: 
\begin{enumerate}
\item In the case of the contrast function $C_1$,
$$c_{ij}=\frac{1}{b(r)}\ind_{(i,j)\in B_r}.$$
\item In the case of the contrast function $C_2$ with two balls of
  radius $r_1<r_2$,
$$c_{ij}=\frac{1}{b(r_1)}\ind_{(i,j)\in B_{r_1}} - 
\frac{1}{b(r_2)}\ind_{(i,j)\in B_{r_2}}.$$
%$$c_{ij}=\begin{cases}
%0 & \mbox{if }(i,j)\not\in B_{r_2}\\
%\frac{1}{b(r_2)} & \mbox{if }(i,j)\in B_{a_2}\setminus B_{r_1}\\
%\frac{1}{b(r_1)}-\frac{1}{b(r_2)} & \mbox{if }(i,j)\in B_{r_1}. 
%\end{cases}$$
\end{enumerate}

Therefore, the computations of the laws will
be similar and they will be treated simultanously using the coefficients
$c_{ij}$.

When the contrast function $C_2$ is used with two radii $r_1 < r_2$, we then  compute the Laplacian $\Delta_r I$ with the larger ball radius, that is with $r=r_2$.

\begin{lemma}
For both contrast functions $C_1$ and $C_2$, the vector
$$\left(C_x(u,v),C_y(u,v),\Delta_r I(u,v)\right)$$
is a Gaussian vector with mean $0$ and covariance matrix of the form:
$$\left( 
\begin{array}{ccc}
\sigma_x^2 & 0 & \sigma_{x,\Delta} \\
0 & \sigma_y^2 & 0 \\
\sigma_{x,\Delta} & 0 & \sigma_{\Delta}^2 
\end{array} \right)$$
In particular, we have that $C_y$ is
independent of $\left(C_x,\Delta_r I\right).$
\end{lemma}

\begin{demo}
The lemma is a consequence of the two following remarks.
The first one is that, in both cases for the contrast function, the coefficients $c_{ij}$ are symmetric: $c_{i,j}=c_{-i,j}$ and $c_{i,j}=c_{i,-j}$.  Thus they satisfy 
that whenever $k$ or $l$ is odd then
\begin{equation}\label{cij=0}
\sum_{(i,j)\in B_r}i^k j^l c_{ij}=0
\end{equation}

The second remark is that the vectors $\vec I_i$ and $\vec I_k$ are independent 
if $i\ne k$. And we thus have
$$\E\left[I_{i,j}I_{k,l}\right]=\begin{cases}
0 & \mbox{if } i \neq k ,\\
\Gamma (j,l) & \mbox{if } i=k . 
\end{cases}$$

We can now compute the covariance matrix. For instance, let us start with:
\begin{align*}
\E\left[C_xC_y\right] & = \sum_{(i,j,k,l)} j k  c_{ij}c_{kl}\E\left[I_{u+i,v+j}I_{u+k,v+l}\right]\\
%& =\sum_{(i,j,l)} i l c_{ij}c_{il}\E\left[I_{u+i,v+j}I_{u+i,v+l}\right]\\
& =\sum_{(i,j,l)} j i c_{ij}c_{il}\Gamma(v+j,v+l)\\
& =\sum_{(j,l)} j \Gamma(v+j,v+l)\sum_{i}i c_{ij}c_{il} = 0 . 
\end{align*}

Similar computations give $\E\left[C_y\Delta_r I\right]  = 0 $ and
\begin{align*}
\sigma_x^2  := & \E\left[C_x^2\right]
 = \sum_{(i,j,l)} j l c_{ij}c_{il}\Gamma(v+j,v+l) ;\\
\sigma_y^2  := & \E\left[C_y^2\right] 
=\sum_{(i,j,l)} i^2 c_{ij}c_{il}\Gamma(v+j,v+l) ; \\
\sigma_\Delta^2 := & \E\left[(\Delta_r I)^2\right] 
=\frac{1}{\beta^2(r)}\sum_{(i,j,l)\in
  \Omega_r}(\alpha(r)+i^2+j^2)(\alpha(r)+i^2+l^2)\Gamma(v+j,v+l) ;\\
\sigma_{x,\Delta} := & \E\left[C_x\Delta_r
  I\right]=\frac{1}{\beta(r_2)}\sum_{(i,j,l)\in
  \Omega_{r_2}}jc_{ij}(\alpha(r)+i^2+l^2)\Gamma(v+j,v+l) ,
\end{align*}
where we have set
$\Omega_r=\{(i,j,l) \text{ such that } (i,j)\in B_r \text{ and } (i,l)\in B_r\}$.

\end{demo}

\subsection{Computation of the threshold}

Now, as we have no more independence between the first and the second order 
derivatives we must compute the conditional law of the contrast function knowing that
 $\Delta_r I=0$.

\begin{proposition}
Let $C$ be one of the two contrast functions. Then, the random variable $\|C\|^2$ is distributed, conditionally on $\{\Delta_r I=0\}$, as the sum of the square of two independent
Gaussian random variables, with mean zero and respective variance
$$\sigma_y^2 \hspace{0.2cm} \text{ and } \hspace{0.2cm} \sigma_{x|\Delta=0}^2=\frac{\sigma_x^2\sigma_\Delta^2-\sigma_{x,\Delta}^2}{\sigma_\Delta^2}, $$
that is a Gamma law with parameters $\frac{1}{2}$ and
$\frac{1}{2}(\sigma_y^2+\sigma_{x|\Delta=0}^2)$.
\end{proposition}

The threshold value $s(\eps)$ defined by
$$\P\bigl(\left\| C \right\| \ge s(\eps) \bigm| \Delta_r I=0\bigr)\le \varepsilon $$
 can no longer be computed explicitly but a
numerical approximation is easy to get as the Gamma density is
well-known.

\begin{demo}
$C_y$ is independent of the pair
$\left(C_x,\Delta_r I\right)$. Thus, conditionally on
$\{\Delta_r I=0\}$, the random variables $C_y$ and $C_x$ are still
independent and the conditional law of $C_y$ is the Gaussian distribution with
mean $0$ and variance $\sigma_y^2$. 

Now, if $D^2:=\sigma_x^2\sigma_\Delta^2-\sigma_{x,\Delta}^2\ne 0$, then
the law of the pair $\left(C_x,\Delta_r I\right)$ has a density which is given by
$$f_{x,\Delta}(t_1,t_2)=\frac{1}{2\pi D}e^{-\frac{1}{2}(t_1,t_2)\Lambda(t_1,t_2)^t}$$
where $\Lambda$ is the inverse of the covariance matrix, i.e.
$$\Lambda=\frac{1}{D^2}\left(\begin{array}{cc}
\sigma_\Delta^2 & -\sigma_{x,\Delta}\\
-\sigma_{x,\Delta} & \sigma_x^2
\end{array}\right).$$
Let us recall that, if $f_\Delta$ denotes the Gaussian density of
$\Delta_r I$, then 
the law of $C_x$ conditionally on $\Delta_r I=0$ has a density given by
$$\frac{f_{x, \Delta}(t_1,0)}{f_\Delta(0)}$$
and so is Gaussian with mean
zero and variance
$$\sigma_{x|\Delta=0}^2=\frac{D^2}{\sigma_\Delta^2}\cdot$$
This result is still valid when $D=0$ since it implies that $C_x$ and $\Delta_r I$ are proportional and thus the law of $C_x$ conditionally on $\Delta_r I=0$ is Gaussian with mean $0$ and variance $0$ (it is not random anymore).
\end{demo}

\subsection{Experiments}

\subsubsection{Case of a piecewise constant object}
\label{expconstant.subsec}

To begin with, we still study the piecewise constant object of Figure
\ref{fig:object} described in Section \ref{sec:experiments_white}. Let
us recall that this image represents a slice of the object that
contains the symmetry axis. The 3-dimensional object is obtained by
rotation around the vertical axis that goes through the middle of the
image.

\begin{figure}[H]
\begin{center}
\begin{tabular}{ccc}
\includegraphics[width=4.2cm]{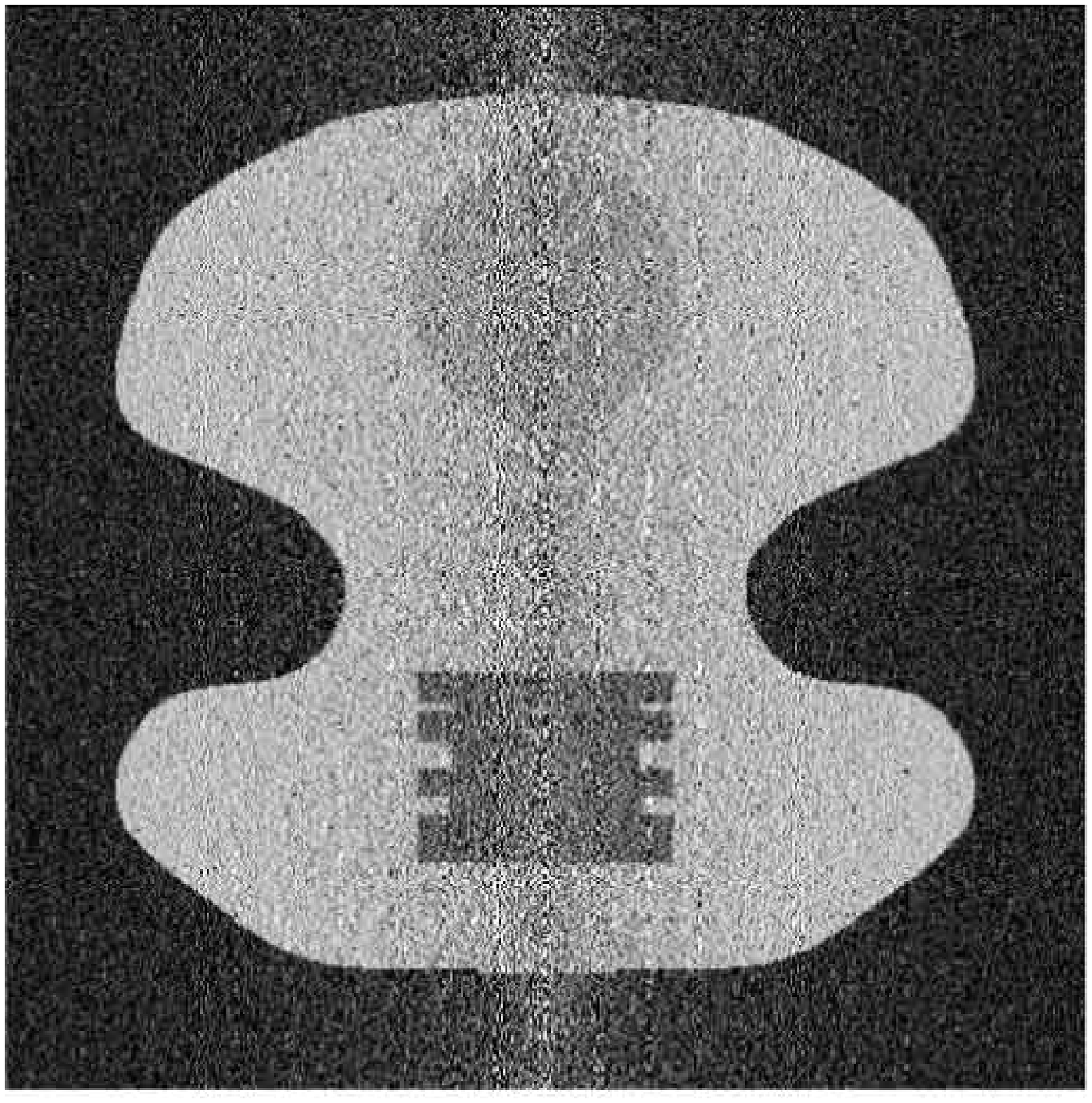} &
\includegraphics[width=4.2cm]{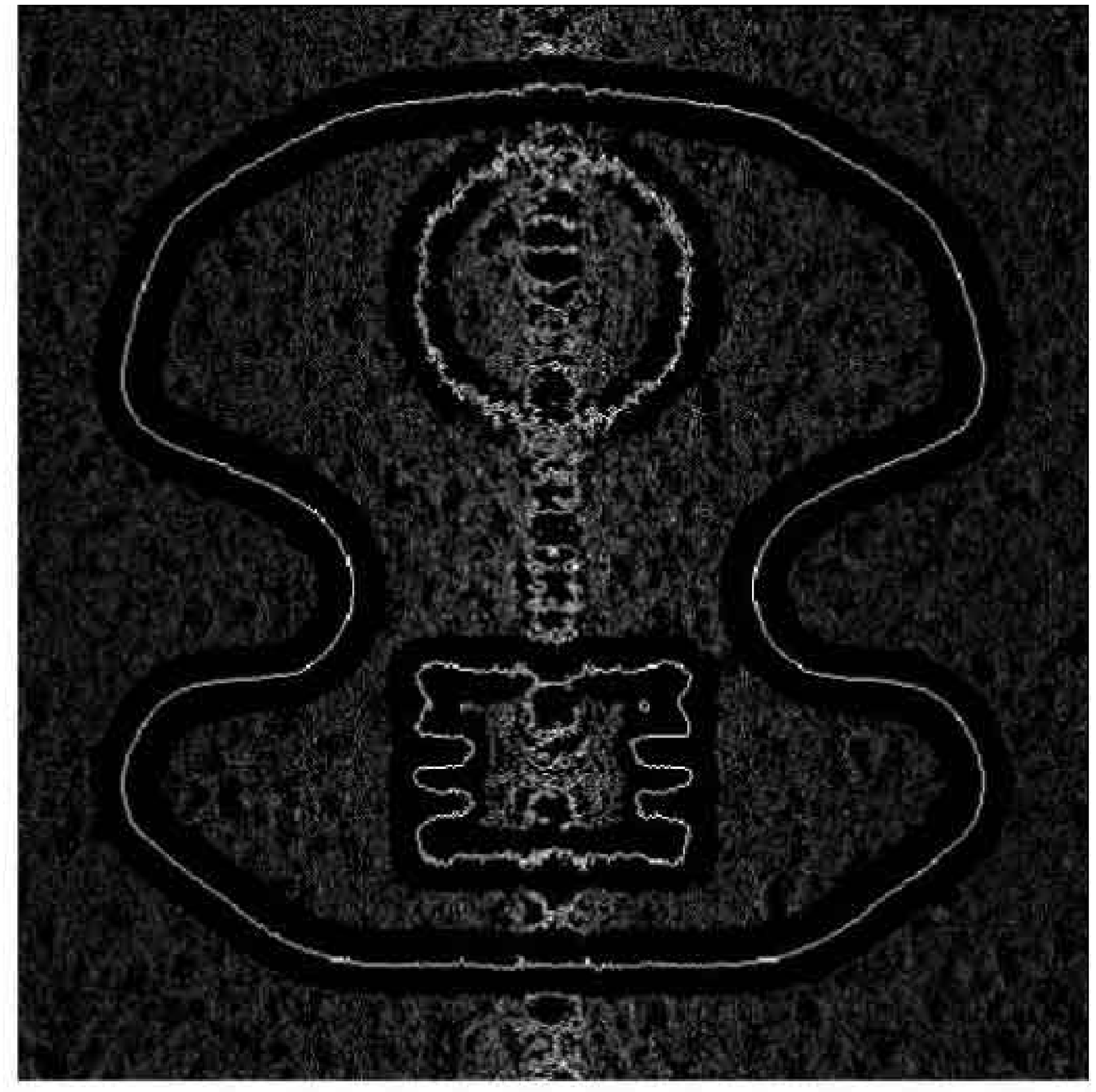} &
\includegraphics[width=4.2cm]{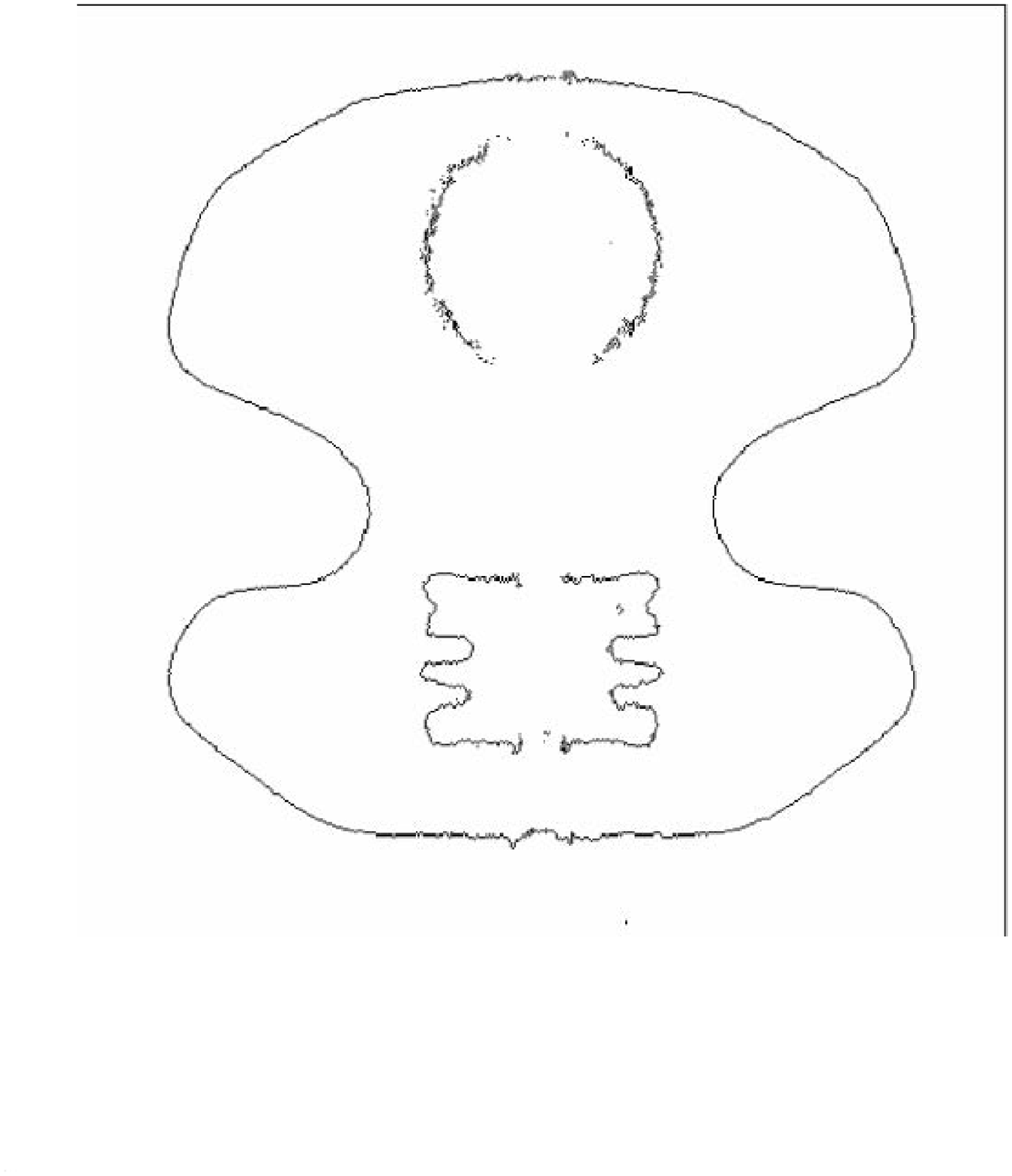}\\
(a) & (b) & (c)
\end{tabular}
\end{center}
\caption{(a) Reconstructed object from a single noisy radiograph, 
(b) Contrast value at the zero-crossings of the Laplacian for the contrast function $C_1$, (c) Significant edges. \label{fig:experiments_tomo_12}}
\end{figure}

In that case, we will use the contrast function $C_1$, which is simply the norm of the gradient. The experiments of Figure  
\ref{fig:experiments_tomo_12} correspond to a ball radius $r=12$ pixels. \\
We start with the image of the radiograph obtained after the
application of matrix $H$ to our initial image. Then a Gaussian white
noise is added to this radiograph. Then tomographic inversion
(application of the matrix $H^{-1}$) is performed. This gives the image of
Figure \ref{fig:experiments_tomo_12}(a). As we already mentioned it, the
noise is not stationary, it is now correlated and its variance
depends on the distance from the symmetry axis. For instance, if the
standard deviation of the Gaussian white noise on the radiograph is
$\sigma=4$, the
variance of the noise on the tomography is about $2\sigma^2=32$ near
the axis, $0.02\sigma^2=0.32$ at a distance of $65$ pixels from the axis
and $8.10^{-3}\sigma^2=0.128$ at the edge of the image located on the right
at $200$ pixels from the axis.

We notice that the edges are not significant near the symmetry axis;
the noise is too important here in order to extract the true edges from
the noise. Let us add that the smaller the difference of the densities
of the material is, the larger the region where the edges are not
significant around the axis is. Even when the edges are significant,
the noise and the method used to detect them can lead to noisy
edges. Moreover, some details are lost because of the smoothing due to
the size of the ball.

\begin{figure}[H]
\begin{center}
\begin{tabular}{ccc}
\includegraphics[width=4.2cm]{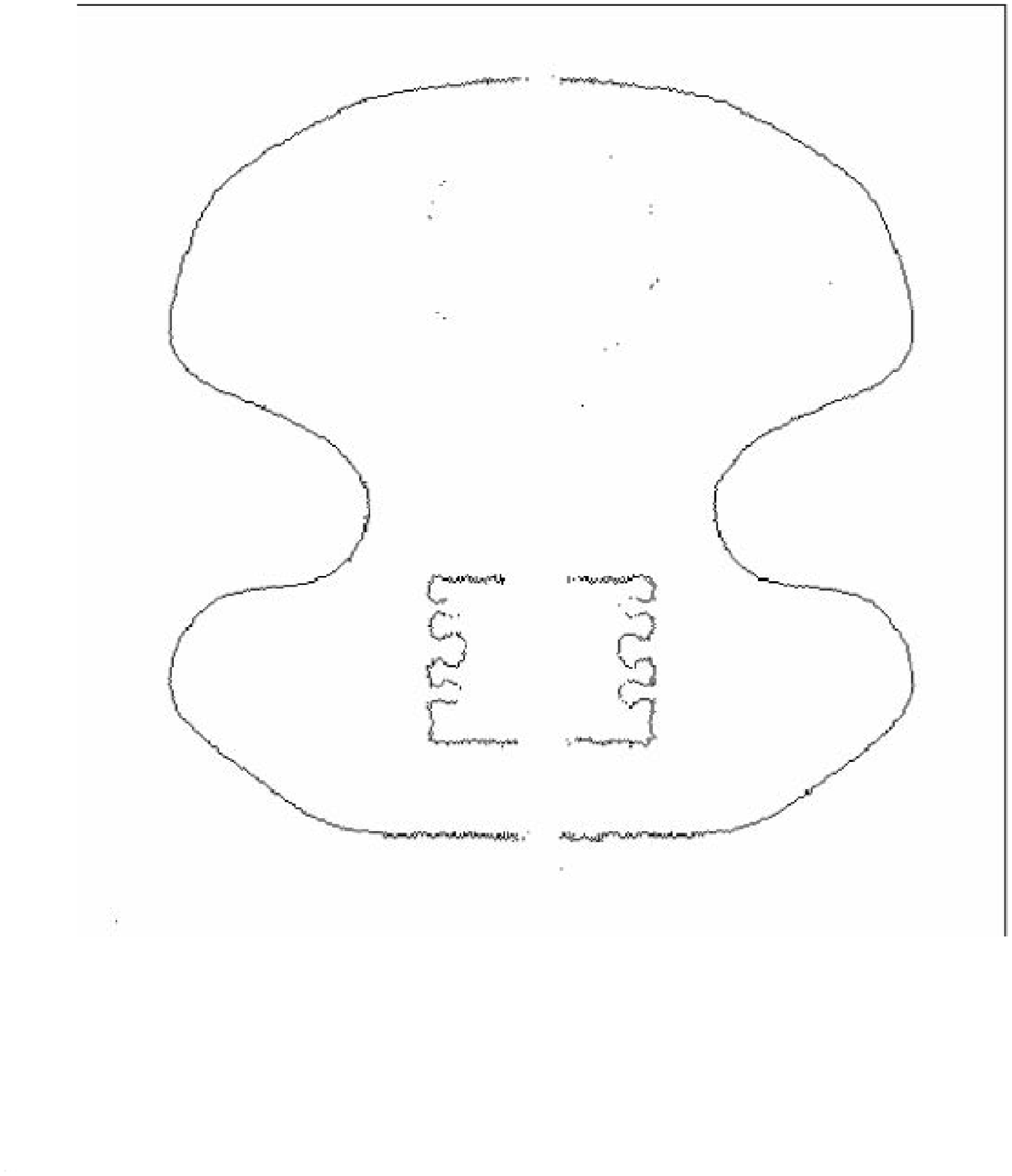} &
\includegraphics[width=4.2cm]{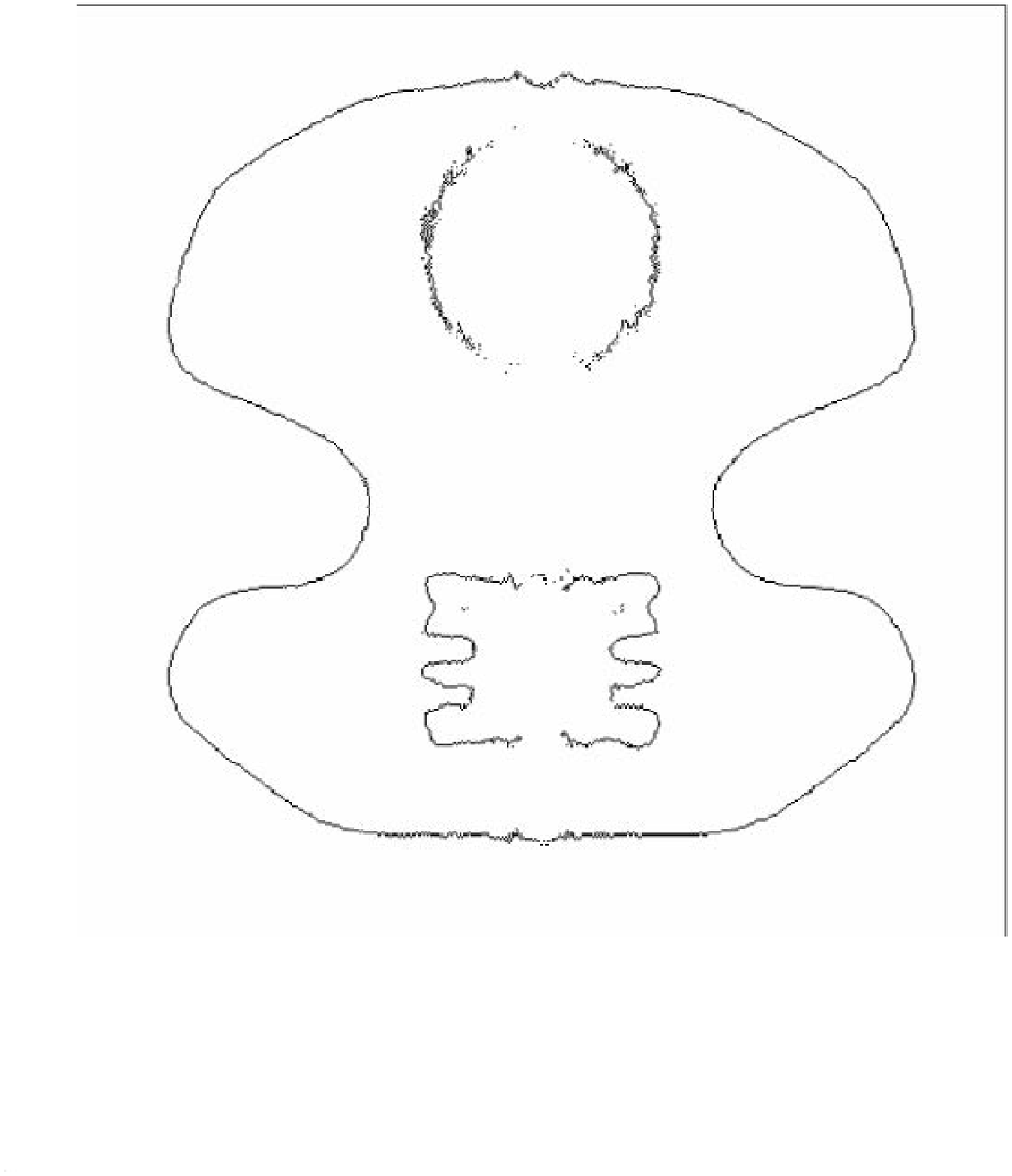} &
\includegraphics[width=4.2cm]{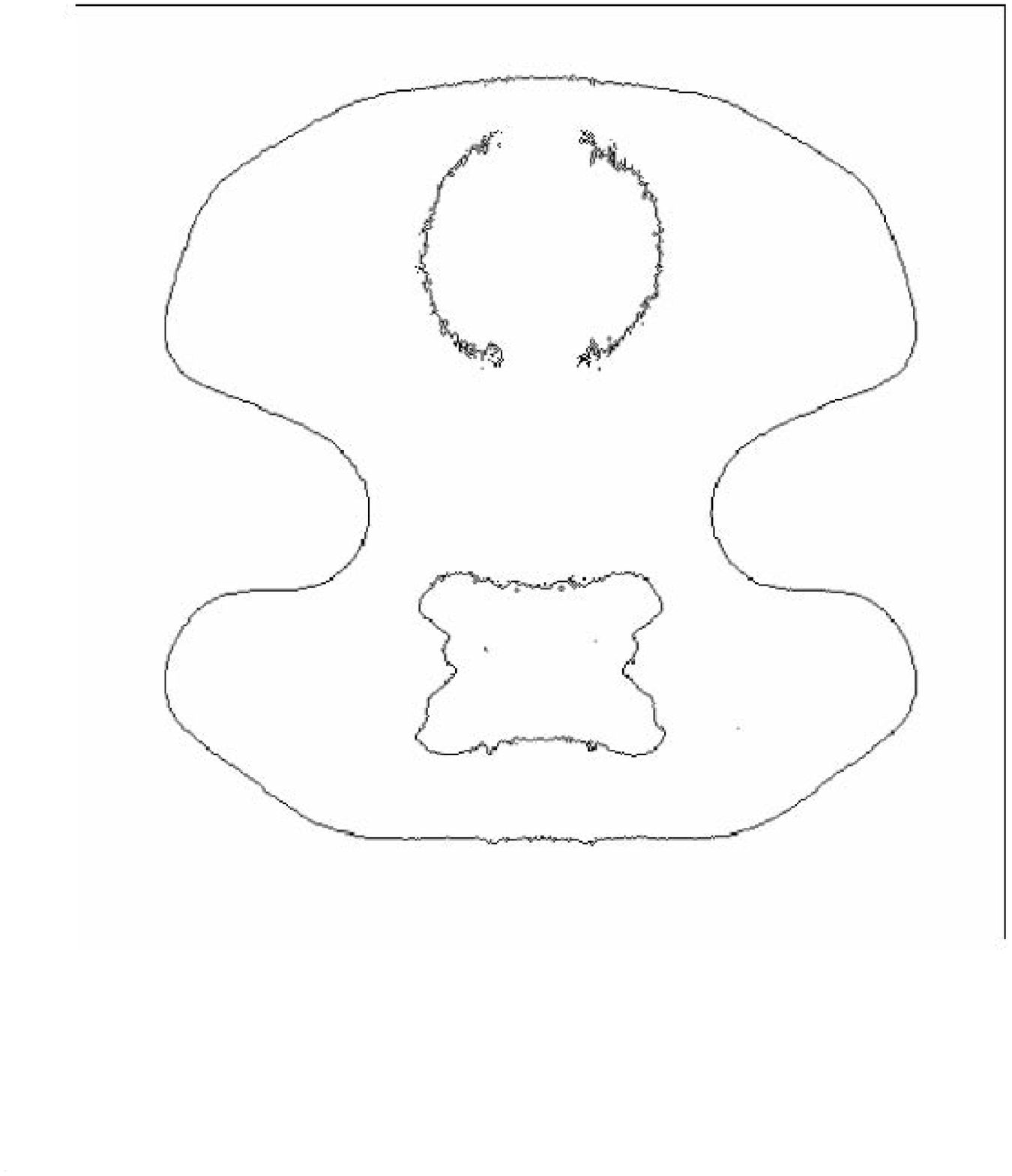}\\
$r=6$ & $r=12$ & $r=20$
\end{tabular}
\end{center}
\caption{\label{fig:different_a}Significant edges obtained with different ball radius: from left to right: $r=6$, $r=12$
and   $r=20$}
\end{figure}

Let us compare the results obtained with different ball radii (see
Figure \ref{fig:different_a}). When the ball radius is small, the edges are more 
accurate but some are not significant: the smoothing of the noise is not
enough to get rid of it. On the contrary, when the radius is large, most of
the edges are 
detected but small details are lost because of this smoothing.

Since the edges separate two materials, one included in another, they
must be closed curves. Usually, an operator has to close them manually.
Our method gives open edges. It does not mean that there is
no edge between the materials: it simply means that the noise level is too high to
give an accurate position of the edge. Therefore, we can then close the
curves manually, or by usual curve completion methods, but this will not tell
which closure is better (i.e. the closest to the real shape).

\begin{figure}[H]
\begin{center}
\begin{tabular}{cc}
\includegraphics[width=5cm]{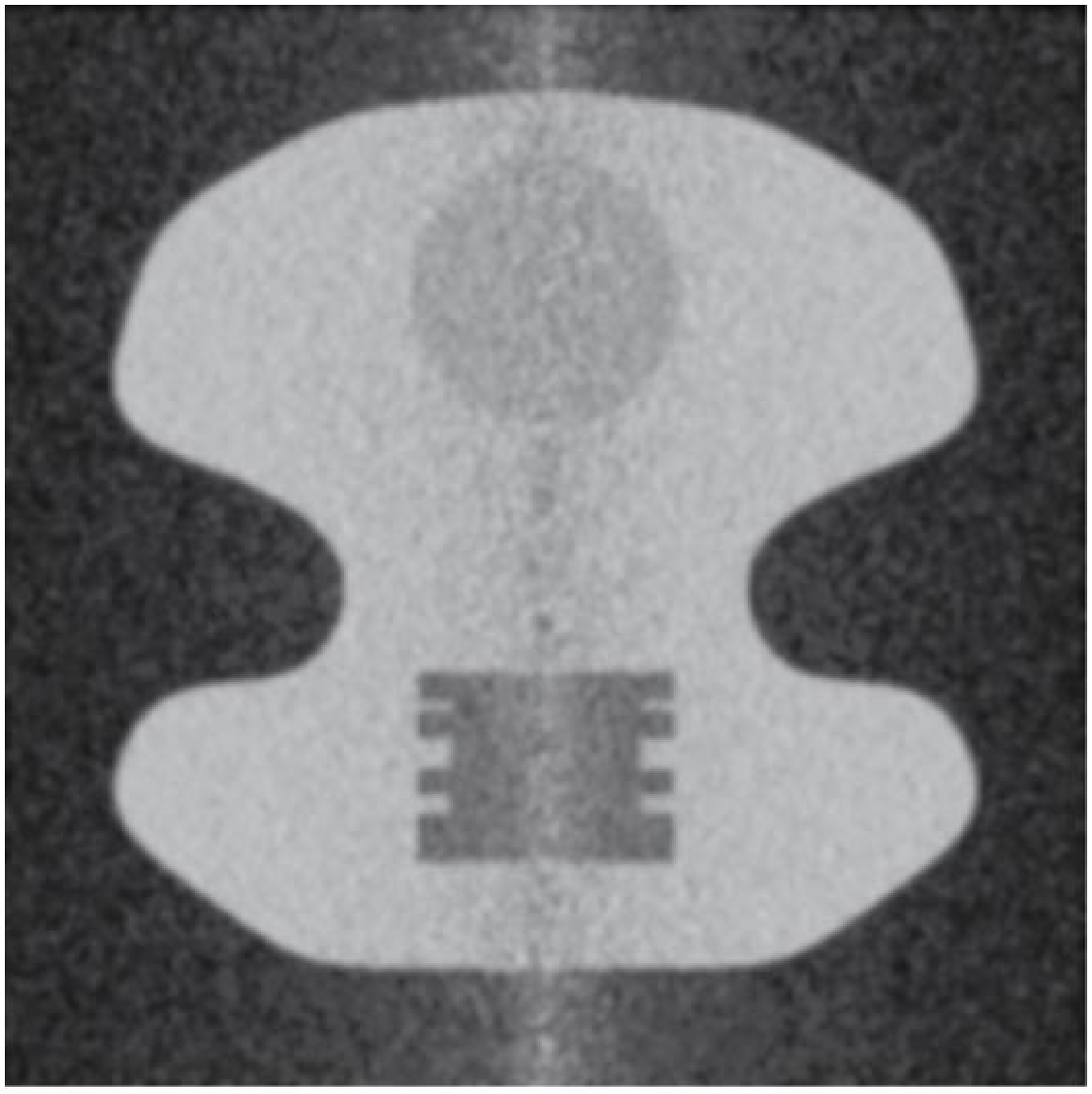} &
\includegraphics[width=5cm]{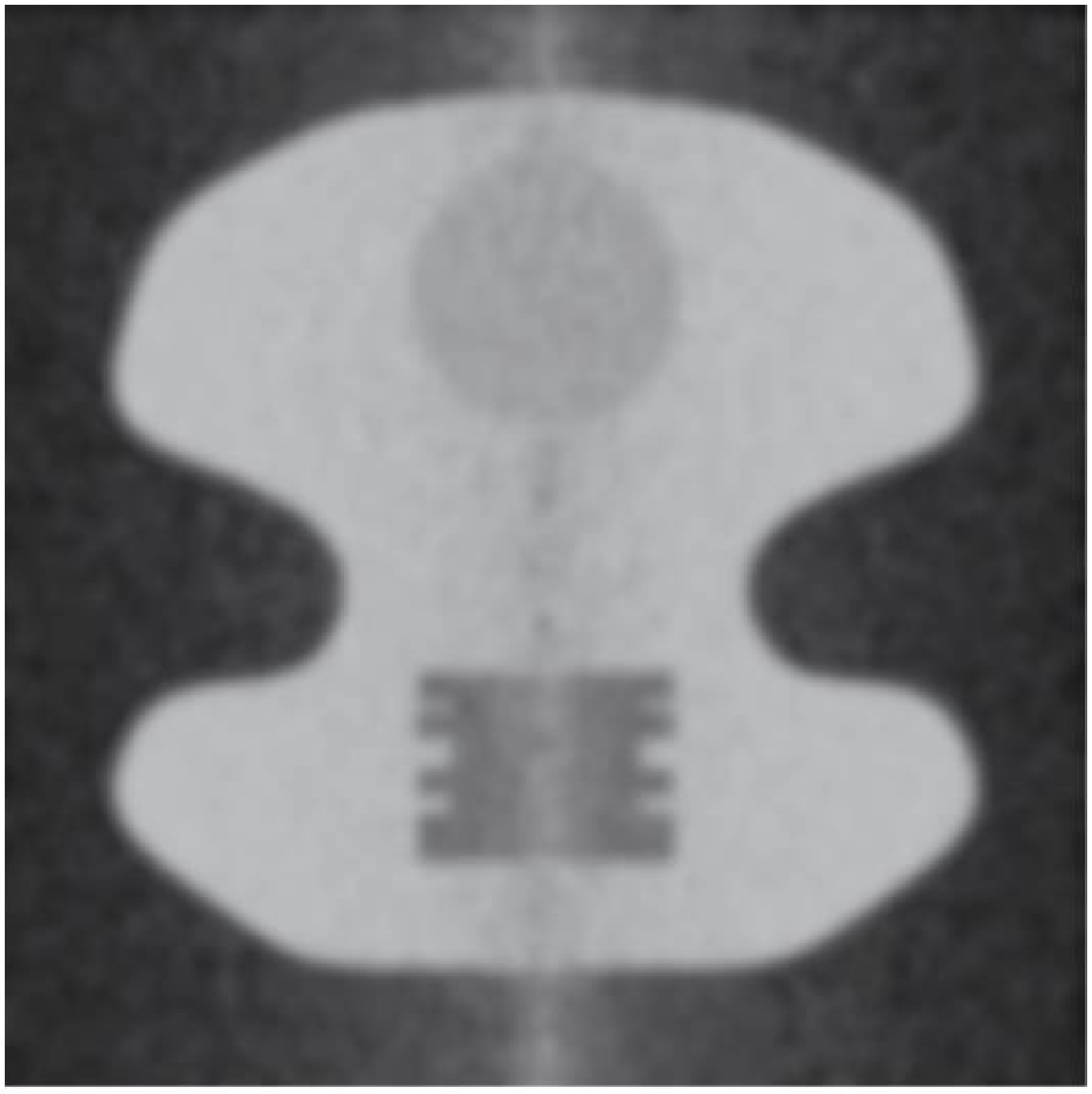} \\
\includegraphics[width=5cm]{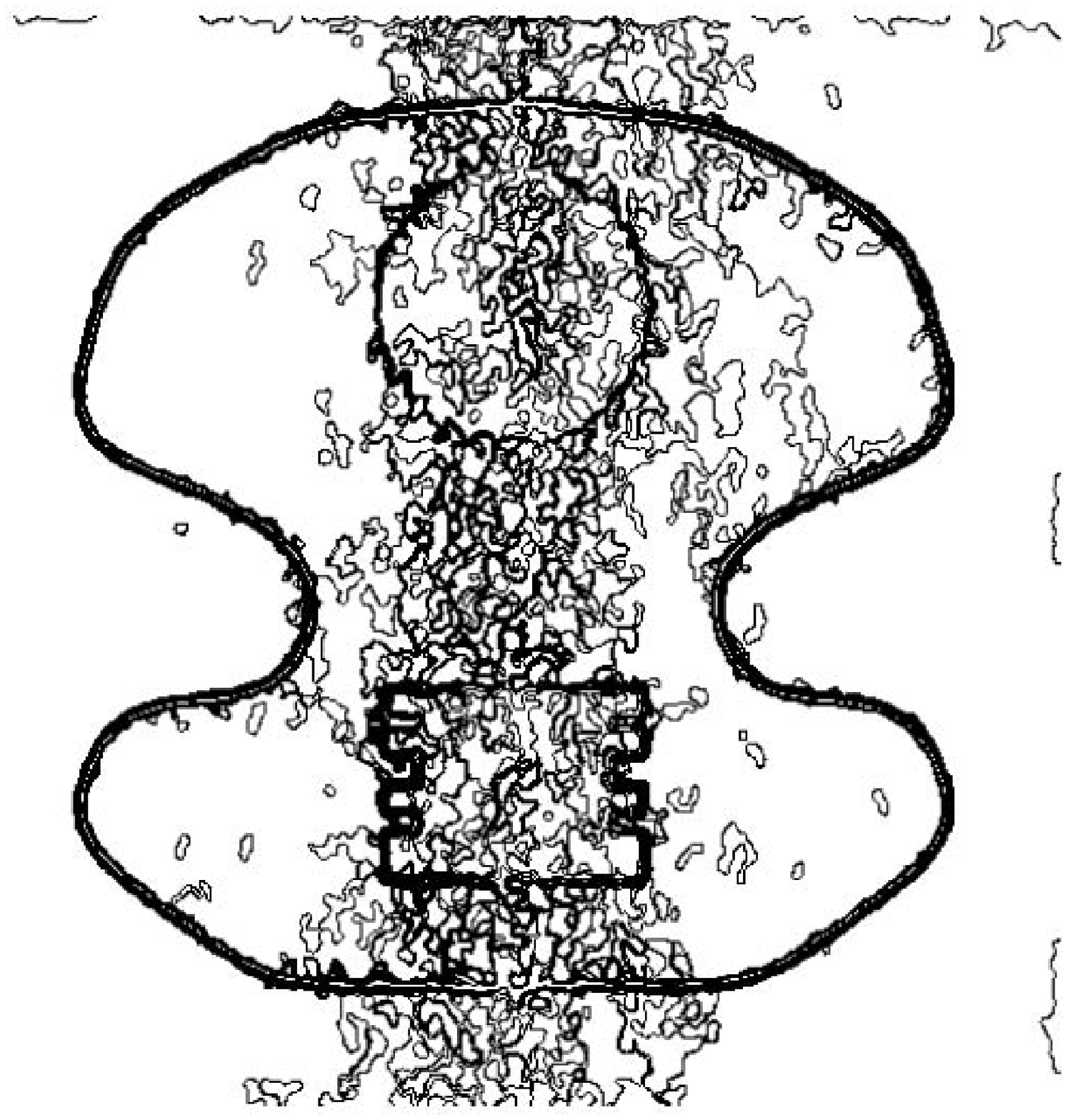} &
\includegraphics[width=5cm]{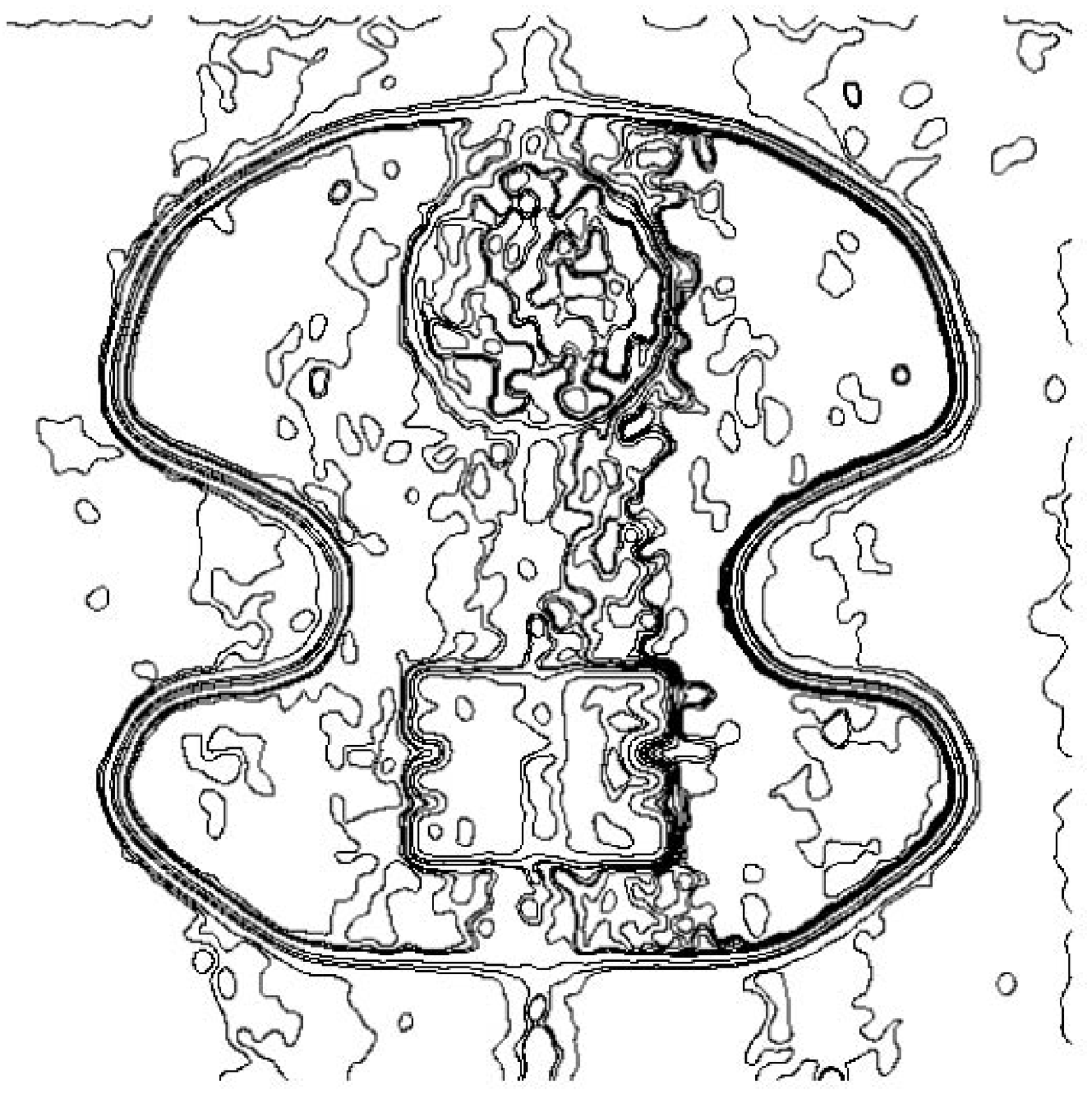}
\end{tabular}
\end{center}
\caption{\label{fig:comparll} First row: the
image of Figure \ref{fig:experiments_tomo_12}(a) is smoothed by convolution with a
Gaussian kernel 
with respective standard deviation $2$ and $4$ pixels. Second
row: the meaningul level lines of each image. }
\end{figure}

\begin{figure}[H]
\begin{center}
\begin{tabular}{ccc}
\includegraphics[width=4.2cm]{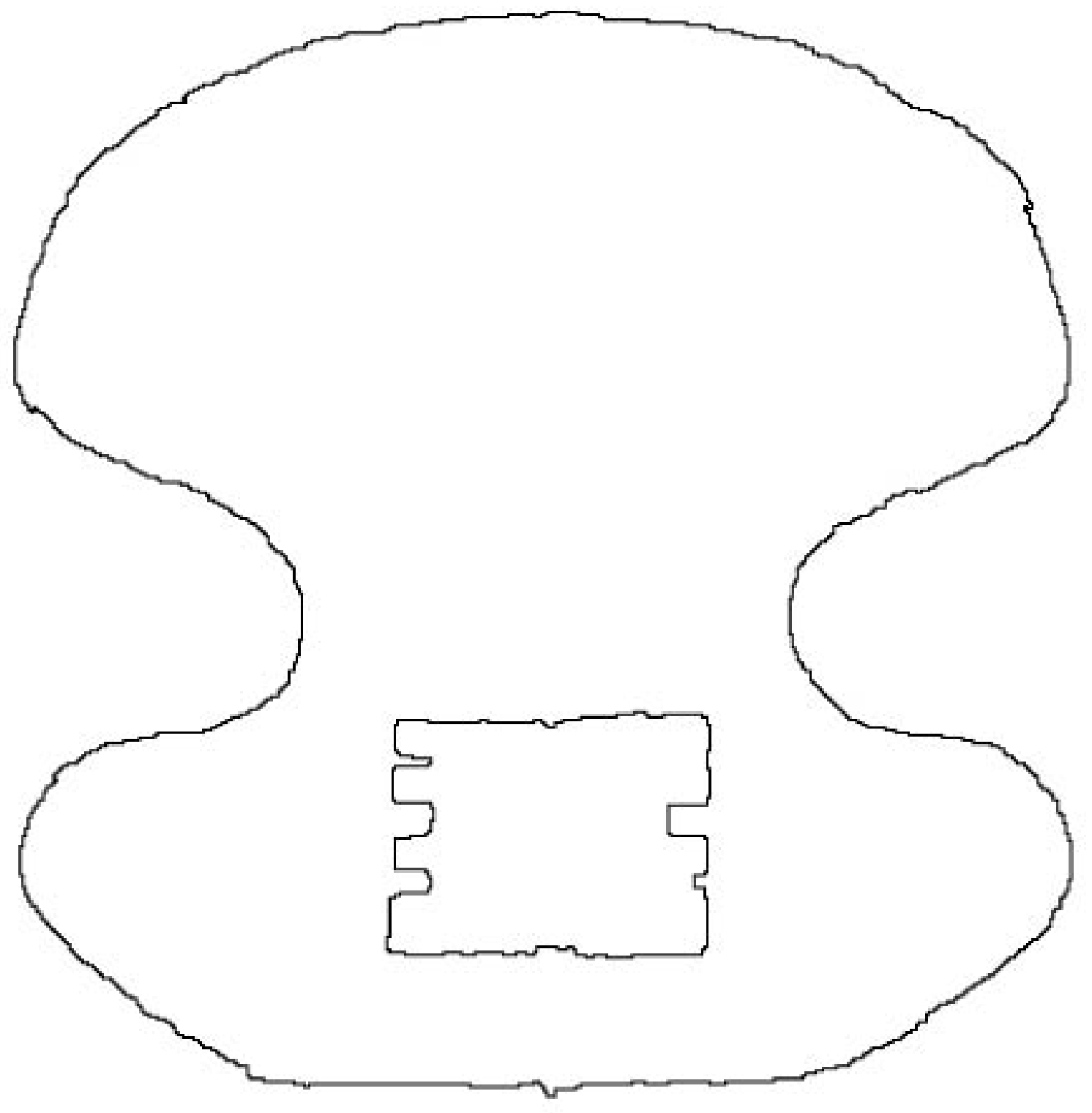} &
\includegraphics[width=4.2cm]{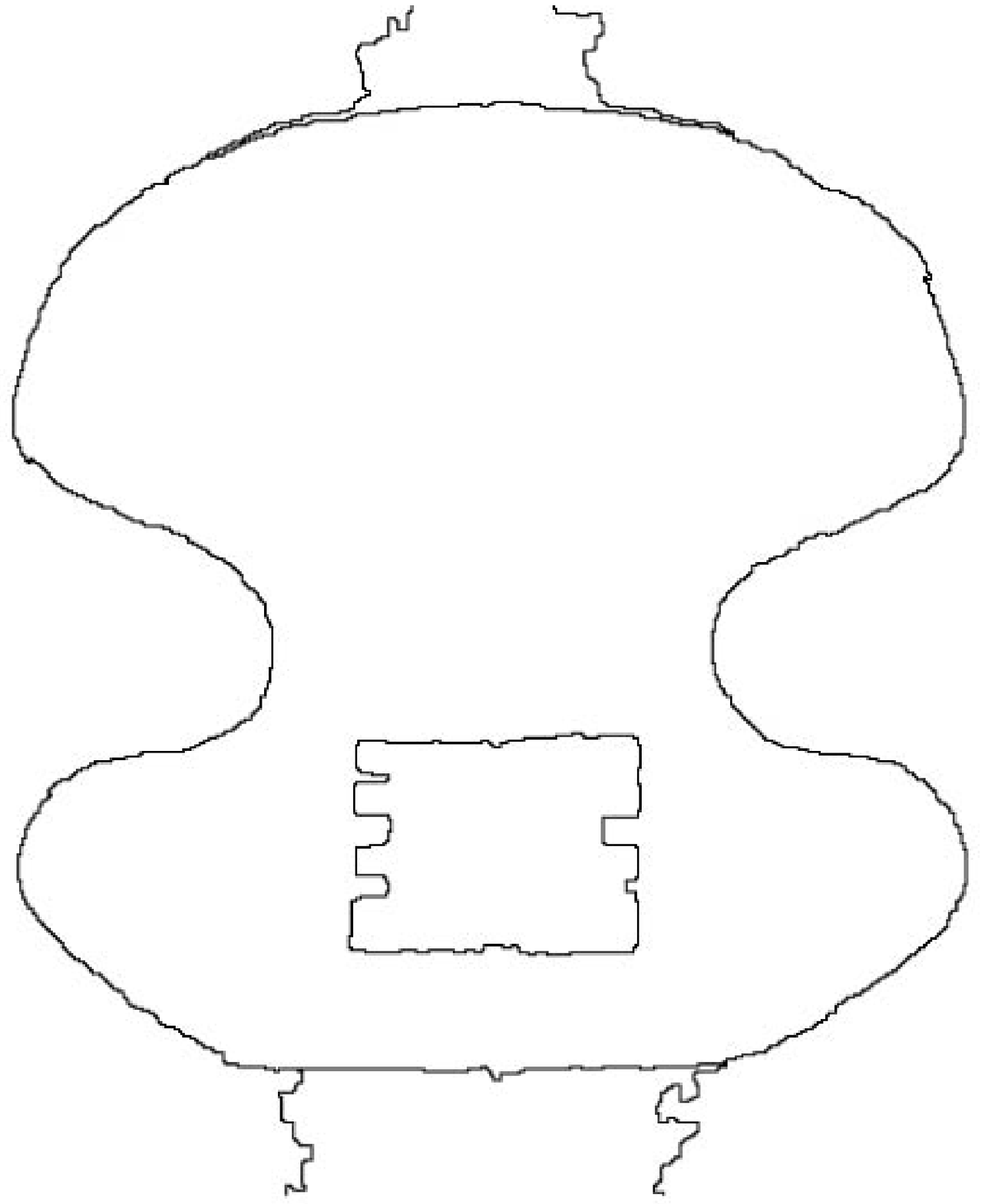} &
\includegraphics[width=4.2cm]{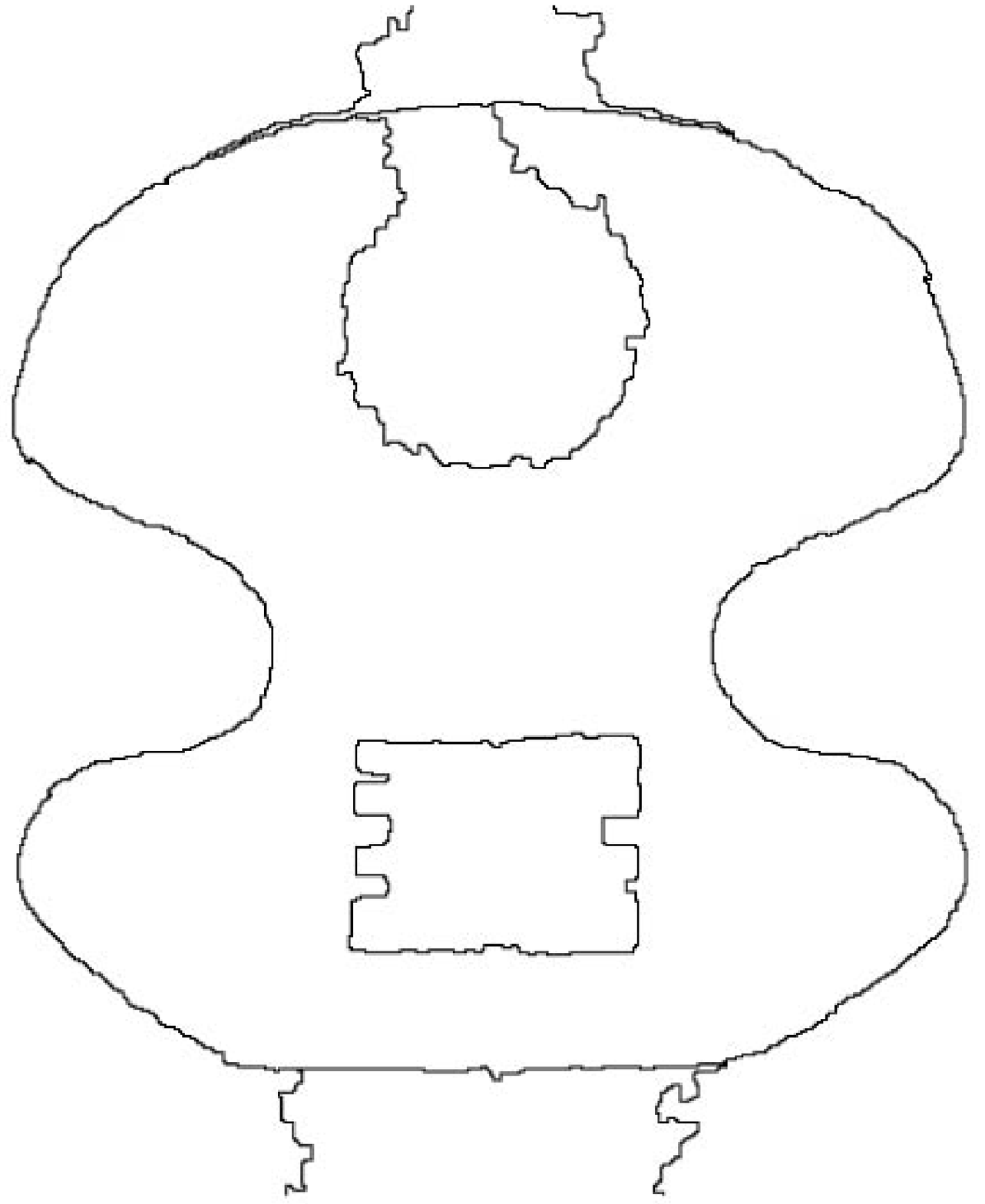} 
\end{tabular}
\end{center}
\caption{\label{fig:comparMS} Results obtained with the Mumford-Shah
  segmentation for piecewise constant images, for three different values of
  $\lambda$. From left to right, the number of regions in the segmented image
  is respectively $3$, $6$ and $7$.  }
\end{figure}

{\it Comparison with other methods.} 
We will give here the results obtain with two other methods which have both
the advantage of directly providing closed curves. \\ 
$\bullet$ The first method is the one introduced in \cite{DMM}. One keeps only the
meaningful level lines of the image, which are defined by: the minimum of the
norm of the gradient along the level line is larger than a threshold
$T(\eps)$. This threshold is computed from the gradient histogram of the
image. The meaning of this definition is that such curves have a probability
less than $\eps$ to appear in a pure noise image (with same gradient
histogram as the original image). The results obtained with this method are
shown on Figure \ref{fig:comparll}. On the first row: we smooth the
image of Figure \ref{fig:experiments_tomo_12}(a) by convolution with a
Gaussian kernel 
with respective standard deviation $2$ and $4$ pixels. And then, on the second
row, we have the respective obtained meaningul level lines. This experiment
clearly shows 
that, since the noise model is not adapted to the image (in particular, the
non-stationarity is not taken into account), many false contours
are detected.\\ 
$\bullet$ The second method is the famous Mumford-Shah segmentation for piecewise
constant images \cite{mumshah1}. Given an observed 
image $g_0$ defined on a domain $D$, one looks for the piecewise
constant approximation $g$ of $g_0$ that minimizes the functional
$$E_{\lambda}(g) = \int_D |g-g_0|^2 + \lambda \, \mathrm{Length}\big(K(g)\big),$$
where $\mathrm{Length}\big(K(g)\big)$ is the one-dimensional 
measure of the discontinuity set of $g$ (which is a set of curves denoted by
$K(g)$) and $\lambda$ is a parameter which weights the second term of the
functional.  The results obtained with this method are shown on Figure
\ref{fig:comparMS} for three different values of $\lambda$. The main drawbacks
of this method are: (a) there is no verification that the obtained contours are
not due to the noise; (b) the parameter $\lambda$ has to be fixed, and the
results are very dependent on its value.

\subsubsection{Case of an inhomogeneous material}

\begin{figure}[H]
\begin{center}
\begin{tabular}{ccc}
\includegraphics[width=4.2cm]{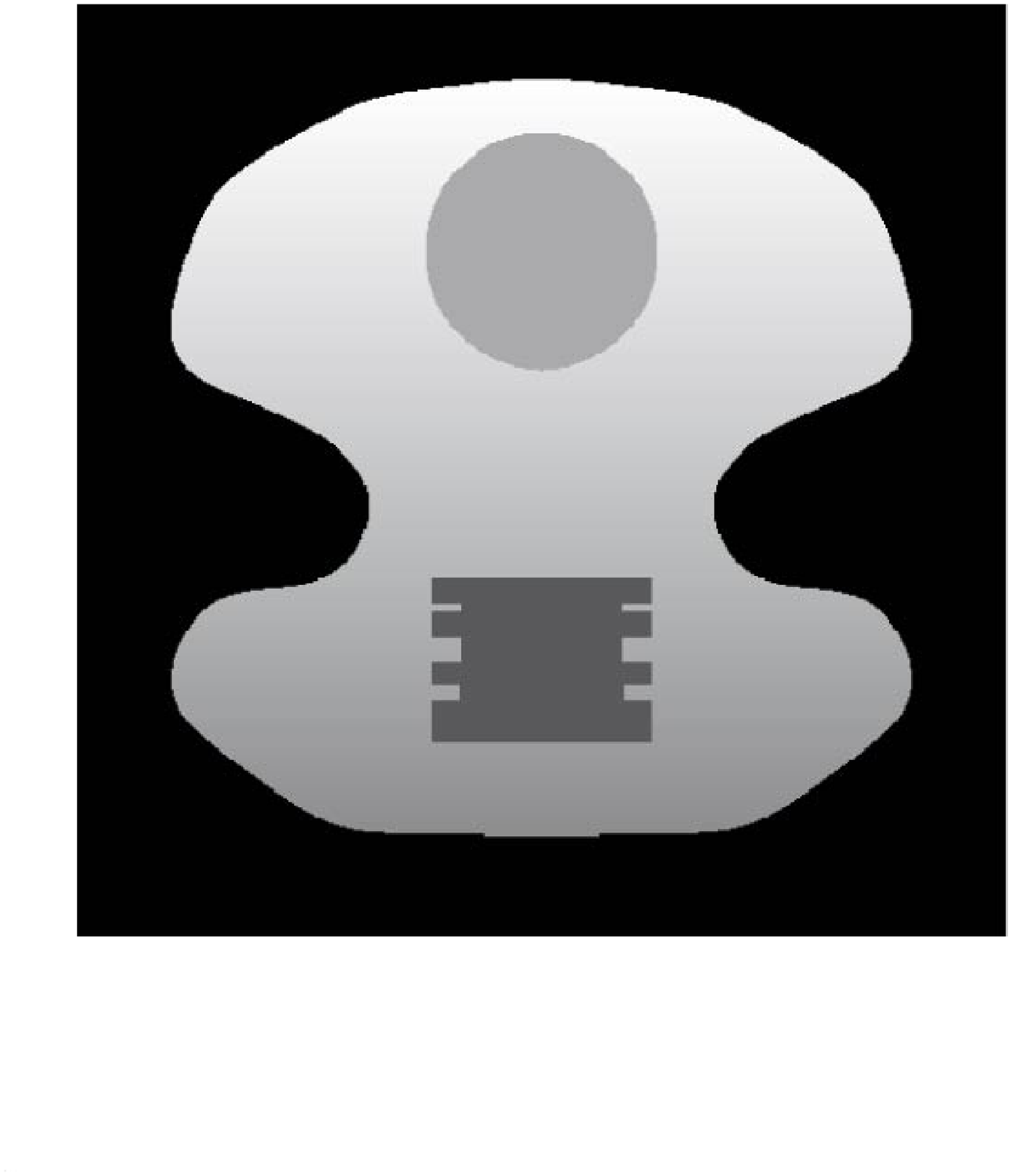} &
\includegraphics[width=4.2cm]{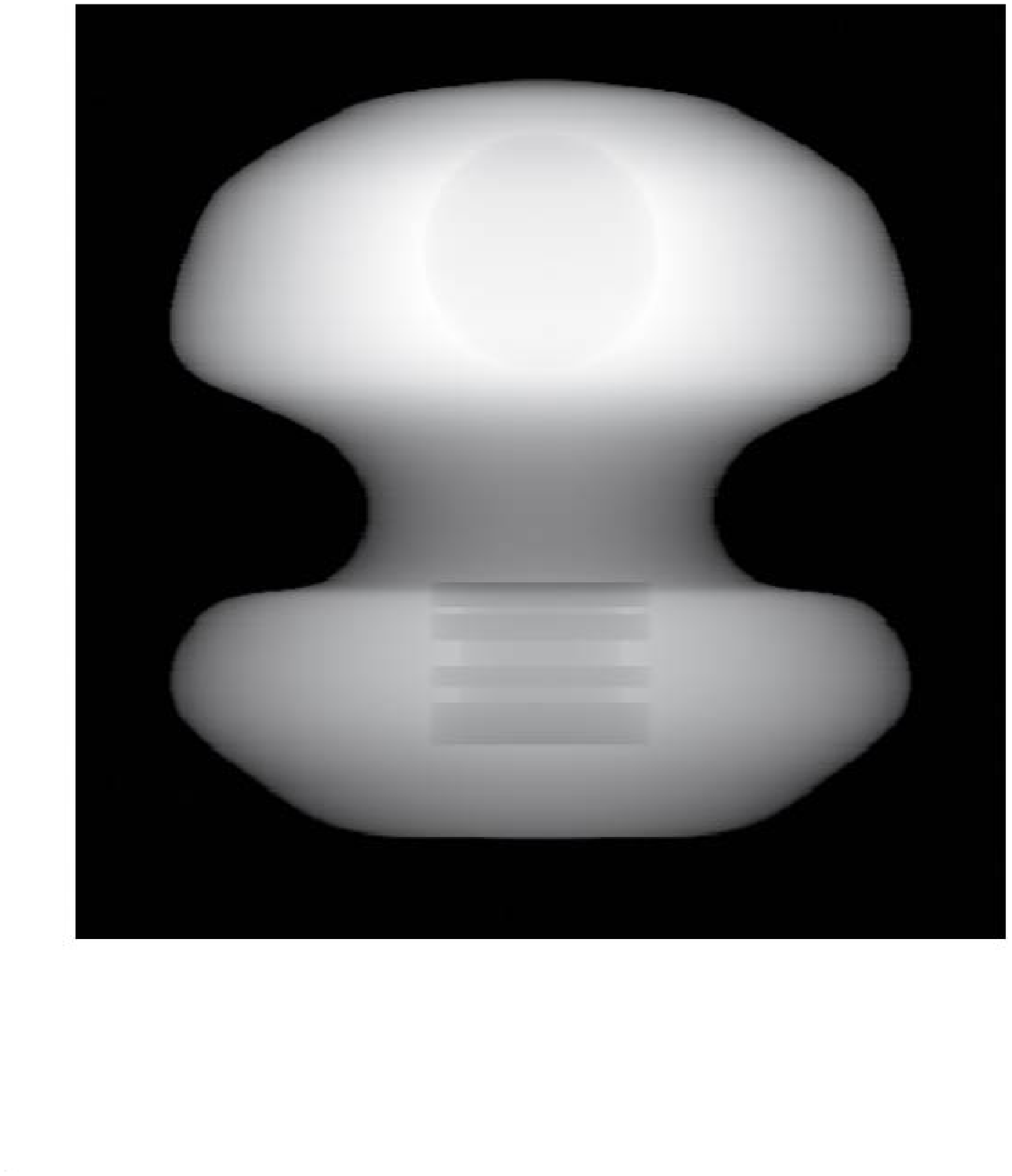} &
\includegraphics[width=4.2cm]{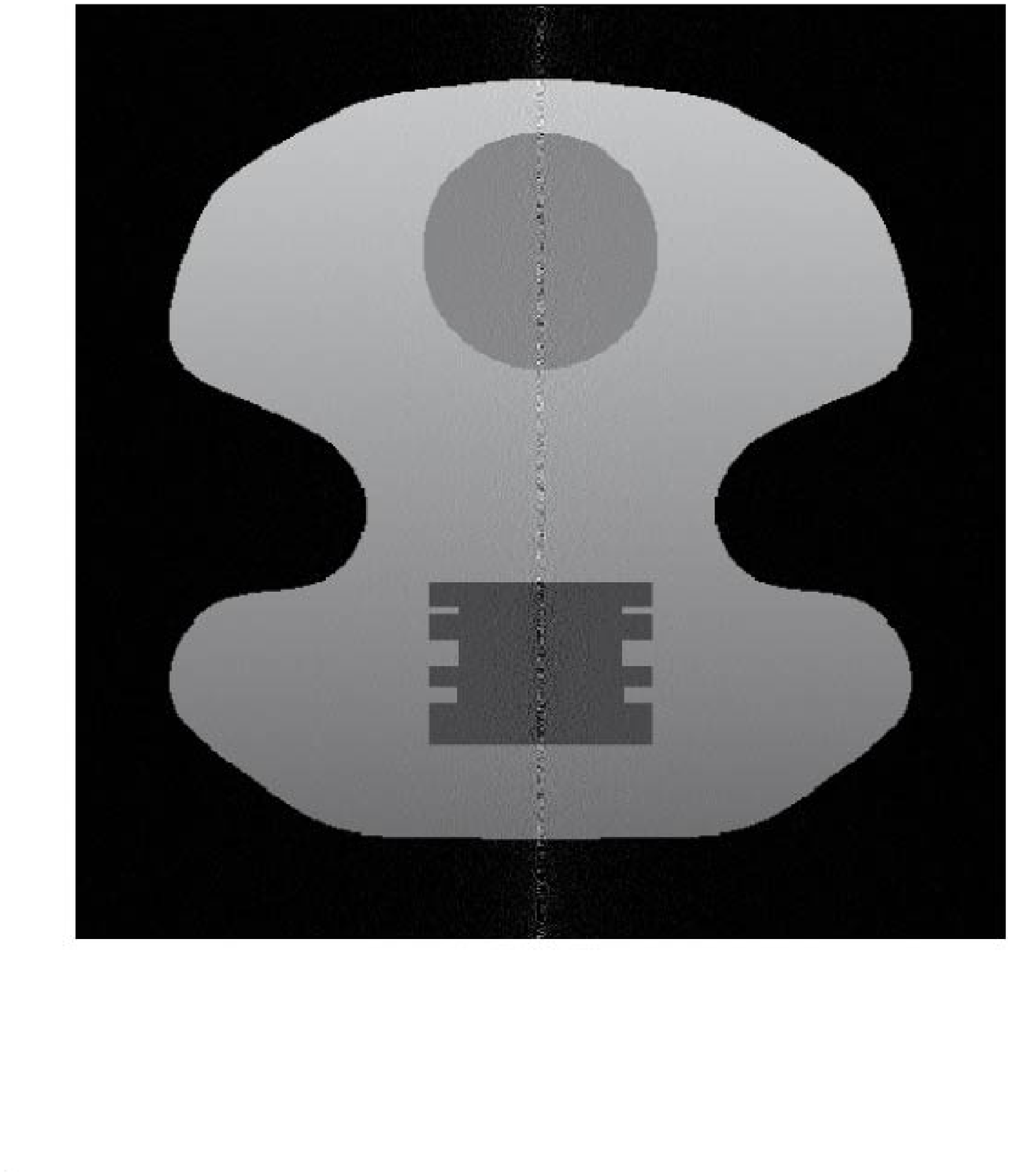}\\
(a) & (b) & (c)
\end{tabular}
\end{center}
\caption{\label{fig:zone}(a) Inhomogeneous object, (b) Its noisy
  radiograph, (c) Tomographic reconstruction}
\end{figure}

\begin{figure}[H]
\begin{center}
\includegraphics[width=5cm]{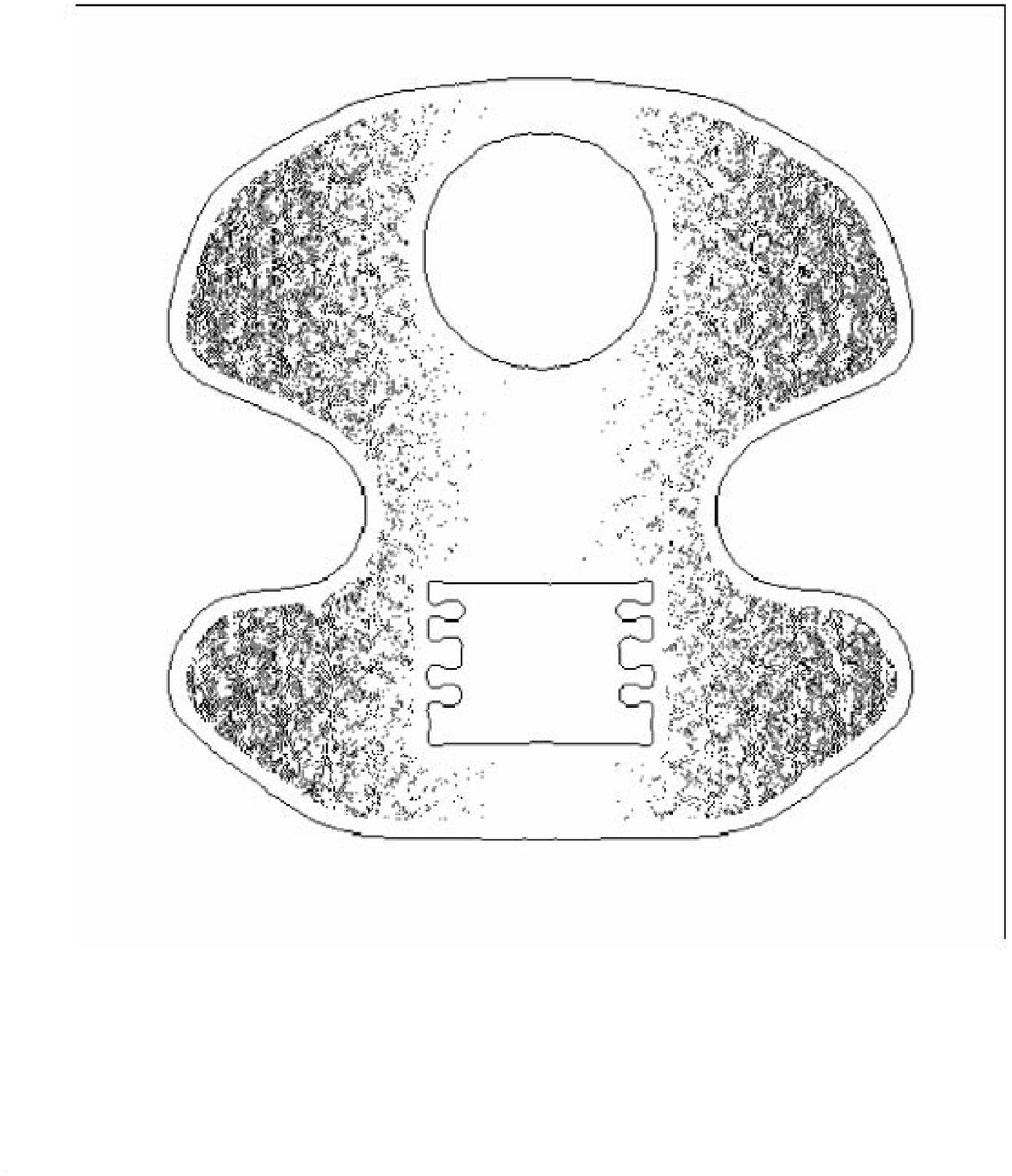}
\end{center}
\caption{\label{fig:C_1+zone}Significant edges with the $C_1$ contrast
  function: there are many false detections.}
\end{figure}

\begin{figure}[H]
\begin{center}
\includegraphics[width=5cm, angle=-90]{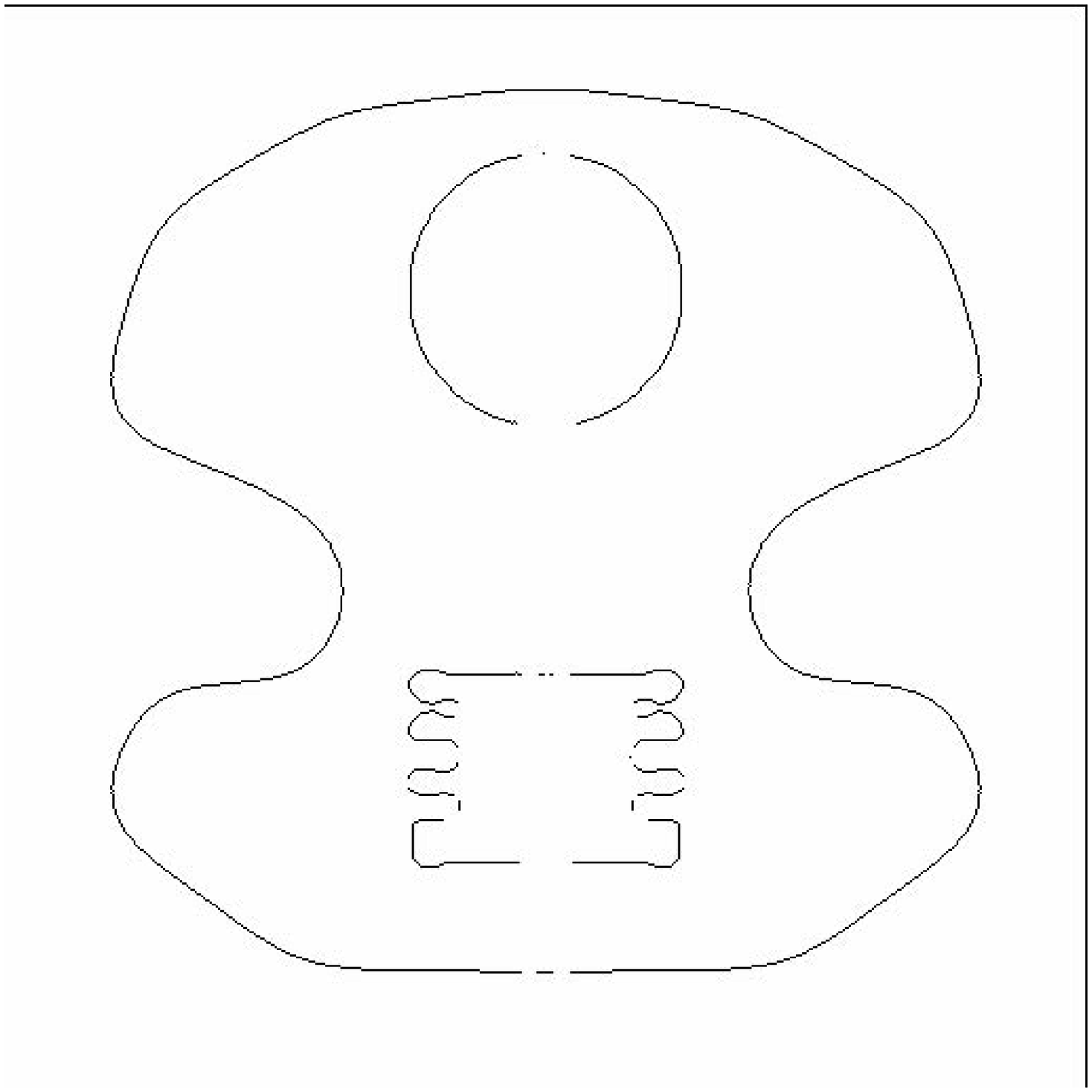}
\end{center}
\caption{\label{fig:C_2+zone}Significant edges with the $C_2$ contrast
  function: only the ``true'' edges are obtained. }
\end{figure}

Let us turn now to a more realistic case: the materials are not
homogeneous and consequently the object is no more piecewise
constant (see Figure \ref{fig:zone}).
As already said, the use of the contrast function $C_1$
fails in that case. This is illustrated by Figure \ref{fig:C_1+zone}.
In this image, one can notice that there are many false detections especially in the
parts of the image where it is not constant.
Figure \ref{fig:C_2+zone} gives the significant edges obtained with
the contrast function $C_2$ with two ball radii $r_1=6$ and
$r_2=12$. With this contrast function, we eventually get only the
``true'' edges.

\begin{thebibliography}{lll}

\bibitem{can}
J. Canny, 
A computational approach to edge detection,
{\it IEEE Trans. on Pattern Analysis and Machine Intelligence} 8, pp. 679-698, 1986.

\bibitem{CMS} 
F. Cao, P. Mus\'e  and F. Sur, 
Extracting Meaningful Curves from Images, 
{\it Journal of Mathematical Imaging and Vision} 22, pp. 159-181, 2005.

\bibitem{Cho} 
I. Abraham, R. Abraham, J.-M. Lagrange and F. Lavallou,
M\'ethodes inverses pour la reconstruction tomographique X monovue,
{\it Revue Chocs} 31 (chocs@cea.fr), 2005.

\bibitem{DMM} 
A. Desolneux, L. Moisan and J.-M. Morel,
 Edge Detection by Helmholtz Principle,
{\it Journal of Mathematical Imaging and Vision} 14, pp. 271-284, 2001.

\bibitem{Din} 
J.M. Dinten,
Tomographie \`a partir d'un nombre limit\'e de projections :
R\'egularisation par des champs markoviens,
PhD Thesis, Universit\'e Paris Sud, 1990.

\bibitem{Hid}
T. Hida,
{\it Brownian Motion},
{Applications of Mathematics} 11,
{Springer-Verlag}, 1980.

\bibitem{HKPS}
T. Hida, H.H. Kuo, J. Potthoff and L. Streit,
{\it White Noise. An infinite Dimensional Calculus.}
{Mathematics and its Applications} 253, 
{Kluwer Academic Publishers Group, Dordrecht}, 1993.

\bibitem{MR} 
D. Marimont and Y. Rubner,
 A probabilistic framework for edge detection and scale selection,
{\it 6th Int. Conference on Computer Vision}, 1998.

\bibitem{mumshah1}
D. Mumford and J. Shah,
Boundary detection by minimizing functionals,
{\it Proc. IEEE Conference on Computer Vision and Pattern Recognition}, San Francisco,
1985.

\bibitem{QB} 
P. Qiu and S. Bhandarkar,
 An edge detection technique using local smoothing and statistical hypothesis testing,
{\it Pattern Recognition Letters} 17, pp. 849-872, 1996.

\bibitem{TLB} 
R. Touzi, A. Lopes and P. Bousquet,
A statistical and geometrical edge detector for SAR images, 
{\it IEEE Transactions on Geoscience and Remote Sensing} 26,
pp. 764-773, 1988.

\bibitem{Wal}
J.B. Walsh,
An introduction to stochastic partial differential equations,
{\it Ecole d'\'et\'e de Probabilit\'es de Saint-Flour XIV 1984}, Lecture
Notes in Math. 1180, {\it Springer, Berlin} 1986.

\end{thebibliography}
\end{document}